\documentclass[12pt]{article}
\usepackage{amsmath}
\usepackage{amsfonts}
\usepackage{amssymb}

\def\thebibliograph#1#2{\section*{{\normalsize \bf #2}}\list
   {[\arabic{enumi}]}{\settowidth\labelwidth{[#1]}\leftmargin\labelwidth
     \advance\leftmargin\labelsep
     \usecounter{enumi}}
     \def\newblock{\hskip .11em plus .33em minus -.07em}
     \sloppy
     \sfcode`\.=1000\relax}

\newtheorem{theorem}{Theorem}
\newtheorem{proposition}{Proposition}
\newtheorem{definition}{Definition}
\newtheorem{corollary}{Corollary}
\newtheorem{lemma}{Lemma}
\newtheorem{remark}{Remark}

\begin{document}

\title{Caffarelli-Kohn-Nirenberg inequalities on Besov and
Triebel-Lizorkin-type spaces}
\author{Douadi Drihem \ \thanks{%
M'sila University, Department of Mathematics, Laboratory of Functional
Analysis and Geometry of Spaces , P.O. Box 166, M'sila 28000, Algeria,
e-mail: \texttt{\ douadidr@yahoo.fr}}}
\date{\today }
\maketitle

\begin{abstract}
We present some Caffarelli-Kohn-Nirenberg-type inequalities on Herz-type
Besov-Triebel-Lizorkin spaces, Besov-Morrey and Triebel-Lizorkin-Morrey
spaces. More precisely, we investigate the inequalities%
\begin{equation*}
\big\|f\big\|_{\dot{k}_{v,\sigma }^{\alpha _{1},r}}\leq c\big\|f\big\|_{\dot{%
K}_{u}^{\alpha _{2},\delta }}^{1-\theta }\big\|f\big\|_{\dot{K}_{p}^{\alpha
_{3},\delta _{1}}A_{\beta }^{s}}^{\theta },
\end{equation*}%
and%
\begin{equation*}
\big\|f\big\|_{\mathcal{E}_{p,2,u}^{\sigma }}\leq c\big\|f\big\|_{\mathcal{M}%
_{\mu }^{\delta }}^{1-\theta }\big\|f\big\|_{\mathcal{N}_{q,\beta
,v}^{s}}^{\theta },
\end{equation*}%
with some appropriate assumptions on the parameters, where $\dot{k}%
_{v,\sigma }^{\alpha _{1},r}$ is the Herz-type Bessel potential spaces,
which are just the Sobolev spaces if $\alpha _{1}=0,1<r=v<\infty $ and $%
\sigma \in \mathbb{N}_{0}$, and $\dot{K}_{p}^{\alpha _{3},\delta
_{1}}A_{\beta }^{s}$ are Besov or Triebel-Lizorkin spaces if $\alpha _{3}=0$
and$\ \delta _{1}=p$. The usual Littlewood-Paley technique, Sobolev and
Franke embeddings are the main tools of this
paper. Some remarks on Hardy-Sobolev inequalities are given.


\end{abstract}

\section{\large Introduction}

Major results in harmonic analysis and partial differential equations invoke
some inequalities. Some examples can be mentioned such as: Caffarelli, Kohn
and Nirenberg in \cite{CKN84}. They proved the following useful inequality:%
\begin{equation}
\big\||x|^{\gamma }f\big\|_{\tau }\leq c\big\||x|^{\beta }f\big\|%
_{q}^{\theta }\big\||x|^{\alpha }\nabla f\big\|_{p}^{1-\theta },\quad f\in 
\mathcal{D}(\mathbb{R}^{n}),  \label{CKN}
\end{equation}%
where $1\leq p,q<\infty ,\tau >0,0\leq \theta \leq 1,\alpha ,\beta ,\gamma
\in \mathbb{R}$ satisfy some suitable conditions. This inequality plays an
important role in theory of PDE's, which extended to fractional Sobolev
spaces by \cite{NM18}. This estimate can be rewritten in the following form:
\begin{equation*}
\big\|f\big\|_{\dot{K}_{\tau }^{\gamma ,\tau }}\leq c\big\|f\big\|_{\dot{K}%
	_{q}^{\beta ,q}}^{\theta }\big\|\nabla f\big\|_{\dot{K}_{p}^{\alpha
		,p}}^{1-\theta },\quad f\in \mathcal{D}(\mathbb{R}^{n}),
\end{equation*}%
where $\dot{K}_{q}^{\alpha ,p}$ is the Herz space, see Definition \ref%
{Herz-spaces} bellow. These function spaces play an important role in
Harmonic Analysis. After they have been introduced in \cite{Herz68}, the
theory of these spaces had a remarkable development in part due to its
usefulness in applications. For instance, they appear in the
characterization of multipliers on Hardy spaces \cite{BS85}, in the  semilinear parabolic equations  \cite{D22}, in the
summability of Fourier transforms \cite{FeichtingerWeisz08}, , and in the Cauchy problem for Navier-Stokes equations \cite{T11}.
For important and latest results on Herz spaces, we refer the reader to the papers \cite{RS}, \cite{ZYZ} and to the monograph \cite{LYH}.

Again \eqref{CKN}\ with $\alpha =\beta =\gamma =0$, is just%
\begin{equation*}
\big\|f\big\|_{M_{\tau }^{\tau }}\leq c\big\|f\big\|_{M_{q}^{q}}^{\theta }\big\|\nabla f\big\|_{M_{p}^{p}}^{1-\theta
},\quad f\in \mathcal{D}(\mathbb{R}^{n}).
\end{equation*}%
where $M_{u}^{p}$, $1\leq u\leq p<\infty $ is the Morrey space.

The main purpose of this paper is to present more general version of such
inequalities. More precisely, we extend this estimate to Herz-type
Besov-Triebel-Lizorkin spaces, called $\dot{K}_{q}^{\alpha ,p}B_{\beta }^{s}$
and $\dot{K}_{q}^{\alpha ,p}F_{\beta }^{s}$, which generalize the usual
Besov and Triebel-Lizorkin spaces. We mean that 
\begin{equation*}
\dot{K}_{p}^{0,p}B_{\beta }^{s}=B_{p,\beta }^{s}\quad \text{and}\quad \dot{K}%
_{p}^{0,p}F_{\beta }^{s}=F_{p,\beta }^{s}.
\end{equation*}%
In addition $\dot{K}_{q}^{\alpha ,p}F_{2}^{0}$ are just the Herz spaces $%
\dot{K}_{q}^{\alpha ,p}$ when $1<p,q<\infty \ $and $-\frac{n}{q}<\alpha <n(1-%
\frac{1}{q})$. In the same manner, we extend these inequalities to
Besov-Morrey and Triebel-Lizorkin-Morrey spaces. Our approach based on the
Littlewood-Paley technique of Triebel \cite{Triebel13} and some results\
obtained by the author in \cite{Drihem1.13,Drihem2.13,Drmana}.

The structure of this paper needs some notation. As usual, $\mathbb{R}^{n}$
denotes the $n$-dimensional real Euclidean space, $\mathbb{N}$ the
collection of all natural numbers and $\mathbb{N}_{0}=\mathbb{N}\cup \{0\}$.
The letter $\mathbb{Z}$ stands for the set of all integer numbers. For any $%
u>0,k\in \mathbb{Z}$ we set $C(u)=\{x\in \mathbb{R}^{n}:\frac{u}{2}%
<\left\vert x\right\vert \leq u\}$ and $C_{k}=C(2^{k})$. $\chi _{k}$, for $%
k\in \mathbb{Z}$, denote the characteristic function of the set $C_{k}$. The
expression $f\approx g$ means that $C$ $g\leq f\leq c\,g$ for some
independent constants $c,C$ and non-negative functions $f$ and $g$.

\noindent For any measurable subset $\Omega \subseteq \mathbb{R}^{n}$ the
Lebesgue space $L^{p}(\Omega )$, $0<p\leq \infty $ consists of all
measurable functions for which 
\begin{equation*}
\big\|f\big\|_{L^{p}(\Omega )}=\Big(\int_{\Omega }\left\vert
f(x)\right\vert ^{p}dx\Big)^{1/p}<\infty ,0<p<\infty
\end{equation*}%
and 
\begin{equation*}
\big\|f \big\|_ {L^{\infty }(\Omega )}=\underset{x\in \Omega }{\text{ess-sup}%
}\left\vert f(x)\right\vert <\infty .
\end{equation*}%
If $\Omega =\mathbb{R}^{n}$, then we put $L^{p}(\mathbb{R}^{n})=L^{p}$ and $%
\big\|f\big\|_{L^{p}(\mathbb{R}^{n})}=\big\|f\big\|_{p}$. The symbol $%
\mathcal{S}(\mathbb{R}^{n})$ is used in place of the set of all Schwartz
functions on $\mathbb{R}^{n}$ and we denote by $\mathcal{S}^{\prime }(\mathbb{R}^{n})$ the dual
space of all tempered distributions on $\mathbb{R}^{n}$. We define the
Fourier transform of a function $f\in \mathcal{S}(\mathbb{R}^{n})$ by 
\begin{equation*}
\mathcal{F(}f)(\xi )=\left( 2\pi \right) ^{-n/2}\int_{\mathbb{R}%
	^{n}}e^{-ix\cdot \xi }f(x)dx, \quad \xi \in \mathbb{R}^{n}.
\end{equation*}%
Its inverse is denoted by $\mathcal{F}^{-1}f$. Both $\mathcal{F}$ and $%
\mathcal{F}^{-1}$ are extended to the dual Schwartz space $\mathcal{S}%
^{\prime }(\mathbb{R}^{n})$ in the usual way. The Hardy-Littlewood maximal
operator $\mathcal{M}$ is defined on $L_{\mathrm{loc}}^{1}$ by%
\begin{equation*}
\mathcal{M}f(x)=\sup_{r>0}\frac{1}{\left\vert B(x,r)\right\vert }%
\int_{B(x,r)}\left\vert f(y)\right\vert dy,\quad x \in \mathbb{R}^{n}
\end{equation*}%
and $\mathcal{M}_{\tau }f=\left( \mathcal{M}\left\vert f\right\vert ^{\tau
}\right) ^{1/\tau }$, $0<\tau <\infty $.

Given two quasi-Banach spaces $X$ and $Y$, we write $X\hookrightarrow Y$ if $%
X\subset Y$ and the natural embedding of $X$ in $Y$ is continuous. We use $c$
as a generic positive constant, i.e.\ a constant whose value may change from
appearance to appearance.

\section{ \large Function spaces}
We start by recalling the definition and some of the properties of the
homogenous Herz spaces $\dot{K}_{q}^{\alpha ,p}$.

\begin{definition}
	\label{Herz-spaces}{\rm Let $\alpha \in \mathbb{R},0<p,q\leq \infty $%
		. The homogeneous Herz space $\dot{K}_{q}^{\alpha ,p}$ is
		defined by %
		\begin{equation*}
		\dot{K}_{q}^{\alpha ,p}=\{f\in L_{\mathrm{loc}}^{q}(\mathbb{R}^{n}\setminus
		\{0\}):\big\|f\big\|_{\dot{K}_{q}^{\alpha ,p}}<\infty \},
		\end{equation*}%
		where %
		\begin{equation*}
		\big\|f\big\|_{\dot{K}_{q}^{\alpha ,p}}=\Big(\sum_{k=-\infty }^{\infty
		}2^{k\alpha p}\big\|f\chi _{k}\big\|_{q}^{p}\Big)^{1/p}
		\end{equation*}%
		with the usual modifications made when $p=\infty $  and/or %
		$q=\infty $.}
\end{definition}

The spaces $\dot{K}_{q}^{\alpha ,p}$ are quasi-Banach spaces and if $\min
(p,q)\geq 1$ then $\dot{K}_{q}^{\alpha ,p}$ are Banach spaces. When $\alpha
=0$ and $0<p=q\leq \infty $ the space $\dot{K}_{p}^{0,p}$ coincides with the
Lebesgue space $L^{p}$. In addition%
\begin{equation*}
\dot{K}_{p}^{\alpha ,p}=L^{p}(\mathbb{R}^{n},|\cdot |^{\alpha p}),\quad 
\text{(Lebesgue space equipped with power weight),}
\end{equation*}%
where%
\begin{equation*}
\big\|f\big\|_{L^{p}(\mathbb{R}^{n},|\cdot |^{\alpha p})}=\Big(\int_{\mathbb{%
		R}^{n}}\left\vert f(x)\right\vert ^{p}|x|^{\alpha p}dx\Big)^{1/p}.
\end{equation*}%
Notice that%
\begin{equation*}
\dot{K}_{q}^{\alpha ,p}\hookrightarrow \mathcal{S}^{\prime }(\mathbb{R}^{n})
\end{equation*}%
for any $\alpha <n(1-\frac{1}{q})$, $1\leq p,q\leq \infty $ or $\alpha = n(1-\frac{1}{q})$, $%
p=1$ and $1\leq q\leq \infty $. We mean that,%
\begin{equation*}
T_{f}(\varphi )=\int_{\mathbb{R}^{n}}f(x)\varphi (x)dx,\quad \varphi \in 
\mathcal{S}(\mathbb{R}^{n}),f\in \dot{K}_{q}^{\alpha ,p}
\end{equation*}%
generates a distribution $T_{f}\in \mathcal{S}^{\prime }(\mathbb{R}^{n})$. A
detailed discussion of the properties of these spaces my be found in \cite%
{HerYang98,LiYang96,LuYangHu}, and references therein.

The following lemma is the $\dot{K}_{q}^{\alpha ,p}$-version of the
Plancherel-Polya-Nikolskij inequality.

\begin{lemma}
	\label{Bernstein-Herz-ine1}\textit{Let }$\alpha _{1},\alpha _{2}\in \mathbb{R%
	}\mathit{\ }$\textit{and} $0<s,\tau ,q,r\leq \infty $. \textit{We suppose
		that }$\alpha _{1}+\frac{n}{s}>0,0<q\leq s\leq \infty $ and $\alpha _{2}\geq
	\alpha _{1}$. \textit{Then there exists a positive constant }$c>0$\textit{\
		independent of }$R$\textit{\ such that for all }$f\in \dot{K}_{q}^{\alpha
		_{2},\delta}\cap \mathcal{S}^{\prime }( \mathbb{R}^{n}) $\textit{\
		with }$\mathrm{supp}$\textit{\ }$\mathcal{F}f\subset \{\xi :|\xi |\leq R\}$%
	\textit{, we have}%
	\begin{equation*}
	\big\|f\big\|_{\dot{K}_{s}^{\alpha _{1},r}}\leq c\text{ }R^{\frac{n}{q}-%
		\frac{n}{s}+\alpha _{2}-\alpha _{1}}\big\|f\big\|_{\dot{K}_{q}^{\alpha
			_{2},\delta }},
	\end{equation*}%
	where%
	\begin{equation*}
	\delta =\left\{ 
	\begin{array}{ccc}
	r, & \text{if} & \alpha _{2}=\alpha _{1}, \\ 
	\tau , & \text{if} & \alpha _{2}>\alpha _{1}.%
	\end{array}%
	\right.
	\end{equation*}
\end{lemma}

\begin{remark}
	{\rm	We would like to mention that Lemma \ref{Bernstein-Herz-ine1}\ improves the
		classical\ Plancherel-Polya-Nikolskij inequality by taking $\alpha
		_{1}=\alpha _{2}=0,r=s$ and using the embedding $\ell^{q}\hookrightarrow
		\ell^{s}$.}
\end{remark}

In the previous lemma we have not treated the case $s<q$. The next lemma
gives a positive answer.

\begin{lemma}
	\label{Bernstein-Herz-ine2}\textit{Let }$\alpha _{1},\alpha _{2}\in \mathbb{R%
	}\mathit{\ }$\textit{and} $0<s,\tau ,q,r\leq \infty $. \textit{We suppose
		that }$\alpha _{1}+\frac{n}{s}>0,0<s\leq q\leq \infty $ and $\alpha _{2}\geq
	\alpha _{1}+\frac{n}{s}-\frac{n}{q}$. \textit{Then there exists a positive
		constant }$c$\textit{\ independent of }$R$\textit{\ such that for all }$f\in 
	\dot{K}_{q}^{\alpha _{2},\delta}\cap \mathcal{S}^{\prime }( \mathbb{R}%
	^{n}) $\textit{\ with }$\mathrm{supp}$\textit{\ }$\mathcal{F}f\subset
	\{\xi :|\xi |\leq R\}$\textit{, we have}%
	\begin{equation*}
	\big\|f\big\|_{\dot{K}_{s}^{\alpha _{1},r}}\leq c\text{ }R^{\frac{n}{q}-%
		\frac{n}{s}+\alpha _{2}-\alpha _{1}}\big\|f\big\|_{\dot{K}_{q}^{\alpha
			_{2},\delta }},
	\end{equation*}%
	where%
	\begin{equation*}
	\delta =\left\{ 
	\begin{array}{ccc}
	r, & \text{if} & \alpha _{2}=\alpha _{1}+\frac{n}{s}-\frac{n}{q}, \\ 
	\tau , & \text{if} & \alpha _{2}>\alpha _{1}+\frac{n}{s}-\frac{n}{q}.%
	\end{array}%
	\right.
	\end{equation*}
\end{lemma}

The proof of these inequalities is given in \cite{Drihem1.13}\textrm{, }%
Lemmas\ 3.10\ and\ 3.14. Let $1<q<\infty $ and $0<p\leq \infty $. If $f$ is
a locally integrable functions on $\mathbb{R}^{n}$ and $-\frac{n}{q}<\alpha
<n(1-\frac{1}{q})$, then%
\begin{equation}
\big\|\mathcal{M}f\big\|_{\dot{K}_{q}^{\alpha ,p}}\leq c\big\|f\big\|_{\dot{K%
	}_{q}^{\alpha ,p}},  \label{convolution estimate}
\end{equation}%
see \cite{LiYang96}.
We need the following lemma, which is basically a consequence of Hardy’s
inequality in the sequence Lebesgue space $\ell^{q}$.

\begin{lemma}
	\label{lem:lq-inequality} Let $0<a<1$ and $0<q\leq \infty $. Let $\left\{
	\varepsilon _{k}\right\} _{k\in \mathbb{N}_{0}}$ be a sequence of positive real
	numbers, such that 
	\begin{equation*}
	\left\Vert \left\{ \varepsilon _{k}\right\} _{k\in \mathbb{N}_{0}}\right\Vert
	_{\ell ^{q}}=I<\infty .
	\end{equation*}%
	Then the sequences $\left\{ \delta _{k}:\delta _{k}=\sum_{j\leq
		k}a^{k-j}\varepsilon _{j}\right\} _{k\in \mathbb{N}_{0}}$ and $\left\{ \eta
	_{k}:\eta _{k}=\sum_{j\geq k}a^{j-k}\varepsilon _{j}\right\} _{k\in \mathbb{N}_{0}}$ belong to $\ell ^{q}$, and 
	\begin{equation*}
	\big\|\left\{ \delta _{k}\right\} _{k\in \mathbb{N}_{0}}\big\|_{\ell ^{q}}+\big\|%
	\left\{ \eta _{k}\right\} _{k\in \mathbb{N}_{0}}\big\|_{\ell ^{q}}\leq c\,I,
	\end{equation*}%
	with $c>0$ only depending on $a$ and $q$.
\end{lemma}
Some of our results of this paper are based on the following result, see Tang and Yang \cite{TY00}.
\begin{lemma}
	\label{Maximal-Inq}Let $1<\beta <\infty ,1<q<\infty $ and $0<p\leq \infty $.
	If $\{f_{j}\}_{j=0}^{\infty }$ is a sequence of locally integrable functions
	on $\mathbb{R}^{n}$ and $-\frac{n}{q}<\alpha <n(1-\frac{1}{q})$, then%
	\begin{equation*}
	\Big\|\Big(\sum_{j=0}^{\infty }(\mathcal{M}f_{j})^{\beta }\Big)^{1/\beta }%
	\Big\|_{\dot{K}_{q}^{\alpha ,p}}\lesssim \Big\|\Big(\sum_{j=0}^{\infty
	}|f_{j}|^{\beta }\Big)^{1/\beta }\Big\|_{\dot{K}_{q}^{\alpha ,p}}.
	\end{equation*}
\end{lemma}

Now, we present the Fourier analytical definition of Herz-type Besov
and Triebel-Lizorkin spaces\ and recall their basic properties. We first
need the concept of a smooth dyadic resolution of unity. Let $\varphi _{0}$\
be a function\ in $\mathcal{S}(\mathbb{R}^{n})$\ satisfying $\varphi
_{0}(x)=1$\ for\ $\left\vert x\right\vert \leq 1$\ and\ $\varphi _{0}(x)=0$\
for\ $\left\vert x\right\vert \geq \frac{3}{2}$.\ We put $\varphi
_{j}(x)=\varphi _{0}(2^{-j}x)-\varphi _{0}(2^{1-j}x)$ for $j=1,2,3,...$.
Then $\{\varphi _{j}\}_{j\in \mathbb{N}_{0}}$\ is a resolution of unity, $%
\sum_{j=0}^{\infty }\varphi _{j}(x)=1$ for all $x\in \mathbb{R}^{n}$.\ Thus
we obtain the Littlewood-Paley decomposition 
\begin{equation*}
f=\sum_{j=0}^{\infty }\mathcal{F}^{-1}\varphi _{j}\ast f
\end{equation*}%
of all $f\in \mathcal{S}^{\prime }(\mathbb{R}^{n})$ $($convergence in $%
\mathcal{S}^{\prime }(\mathbb{R}^{n}))$.

We are now in a position to state the definition of Herz-type Besov and
Triebel-Lizorkin spaces.

\begin{definition}
	\label{Herz-Besov-Triebel}{\rm Let $\alpha ,s\in \mathbb{R},0<p,q\leq \infty $\
		and $0<\beta \leq \infty $. \newline
		$\mathrm{(i)}$ The Herz-type Besov space $\dot{K}_{q}^{\alpha ,p}B_{\beta
		}^{s}$\ is the collection of all $f\in \mathcal{S}^{\prime }(\mathbb{R}^{n})$%
		\ such that 
		\begin{equation*}
		\big\|f\big\|_{\dot{K}_{q}^{\alpha ,p}B_{\beta }^{s}}=\Big(%
		\sum\limits_{j=0}^{\infty }2^{js\beta }\big\|\mathcal{F}^{-1}\varphi
		_{j}\ast f\big\|_{\dot{K}_{q}^{\alpha ,p}}^{\beta }\Big)^{1/\beta }<\infty ,
		\end{equation*}%
		with the obvious modification if\textit{\ }$\beta =\infty .$\newline
		$\mathrm{(ii)}$ Let $0<p,q<\infty $. The Herz-type Triebel-Lizorkin space $%
		\dot{K}_{q}^{\alpha ,p}F_{\beta }^{s}$ is the collection of all $f\in 
		\mathcal{S}^{\prime }(\mathbb{R}^{n})$\ such that%
		\begin{equation*}
		\big\|f\big\|_{\dot{K}_{q}^{\alpha ,p}F_{\beta }^{s}}=\Big\|\Big(%
		\sum\limits_{j=0}^{\infty }2^{js\beta }\left\vert \mathcal{F}^{-1}\varphi
		_{j}\ast f\right\vert ^{\beta }\Big)^{1/\beta }\Big\|_{\dot{K}_{q}^{\alpha
				,p}}<\infty ,
		\end{equation*}%
		with the obvious modification if\textit{\ }$\beta =\infty .$}
\end{definition}
\begin{remark}
	{\rm	Let\ $s\in \mathbb{R},0<p,q\leq \infty ,0<\beta \leq \infty $ and $\alpha >-%
		\frac{n}{q}$. The spaces\ $\dot{K}_{q}^{\alpha ,p}B_{\beta }^{s}$ and $\dot{K%
		}_{q}^{\alpha ,p}F_{\beta }^{s}$\ are independent of the particular choice
		of the smooth dyadic resolution of unity\ $\{\varphi _{j}\}_{j\in \mathbb{N}%
			_{0}}$ $($in the sense of\ equivalent quasi-norms$)$. In particular $\dot{K}%
		_{q}^{\alpha ,p}B_{\beta }^{s}$ and $\dot{K}_{q}^{\alpha ,p}F_{\beta }^{s}$\
		are quasi-Banach spaces and if\ $p,q,\beta \geq 1$, then they\ are Banach
		spaces. Further results, concerning, for instance, lifting properties,
		Fourier multiplier and local means characterizations can be found in
		\cite{DjDr18, Drihem1.13,Drihem2.13,Drmana,drihem2016jawerth, XuYang03,XuYang05, Xu13}.}
\end{remark}

Now we give the definitions of the spaces $B_{p,\beta }^{s}$ and $F_{p,\beta
}^{s}$.

\begin{definition}
	{\rm	$\mathrm{(i)}$\textit{\ }Let $s\in \mathbb{R}$\ and $0<p,\beta \leq \infty $%
		. The Besov space $B_{p,\beta }^{s}$ is the collection of all $f\in \mathcal{%
			S}^{\prime }(\mathbb{R}^{n})$\ such that 
		\begin{equation*}
		\big\|f\big\|_{B_{p,\beta }^{s}}=\Big(\sum\limits_{j=0}^{\infty }2^{js\beta }%
		\big\|\mathcal{F}^{-1}\varphi _{j}\ast f\big\|_{p}^{\beta }\Big)^{1/\beta
		}<\infty ,
		\end{equation*}%
		with the obvious modification if\textit{\ }$\beta =\infty .\newline
		\mathrm{(ii)}$\textit{\ }Let\textit{\ }$s\in \mathbb{R},0<p<\infty $\ and $%
		0<\beta \leq \infty $\textit{. }The Triebel-Lizorkin space\textit{\ }$%
		F_{p,\beta }^{s}$ is the collection of all $f\in \mathcal{S}^{\prime }(%
		\mathbb{R}^{n})$\textit{\ }such that%
		\begin{equation*}
		\big\|f\big\|_{F_{p,\beta }^{s}}=\Big\|\Big(\sum\limits_{j=0}^{\infty
		}2^{js\beta }\left\vert \mathcal{F}^{-1}\varphi _{j}\ast f\right\vert
		^{\beta }\Big)^{1/\beta }\Big\|_{p}<\infty ,
		\end{equation*}%
		with the obvious modification if\textit{\ }$\beta =\infty .$}
\end{definition}

The theory of the spaces $B_{p,\beta }^{s}$ and $F_{p,\beta }^{s}$ has been
developed in detail in \cite{Triebel83,Triebel92} but has a
longer history already including many contributors; we do not want to
discuss this here. Clearly, for $s\in \mathbb{R},0<p<\infty $ and $0<\beta
\leq \infty ,$%
\begin{equation*}
\dot{K}_{p}^{0,p}B_{\beta }^{s}=B_{p,\beta }^{s}\quad \text{and}\quad \dot{K}%
_{p}^{0,p}F_{\beta }^{s}=F_{p,\beta }^{s}.
\end{equation*}

Let $w$ denote a positive, locally integrable function and $0<p<\infty $.
Then the weighted Lebesgue space $L^{p}(\mathbb{R}^{n},w)$ contains all
measurable functions such that 
\begin{equation*}
\big\|f\big\|_{L^{p}(\mathbb{R}^{n},w)}=\Big(\int_{\mathbb{R}^{n}}\left\vert
f(x)\right\vert ^{p}w(x)dx\Big)^{1/p}<\infty .
\end{equation*}%
For $\varrho \in \lbrack 1,\infty )$ we denote by $\mathcal{A}_{\varrho }$
the Muckenhoupt class of weights, and $\mathcal{A}_{\infty }=\cup _{\varrho
	\geq 1}\mathcal{A}_{\varrho }$. We refer to \cite{GR85} for the general
properties of these classes. Let $w\in \mathcal{A}_{\infty }$, $s\in \mathbb{%
	R}$, $0<\beta \leq \infty $ and $0<p<\infty $. We define weighted Besov
spaces $B_{p,\beta }^{s}(\mathbb{R}^{n},w)$ to be the set of all
distributions $f\in \mathcal{S}^{\prime }(\mathbb{R}^{n})$ such that%
\begin{equation*}
\big\|f\big\|_{B_{p,\beta }^{s}(\mathbb{R}^{n},w)}=\Big(\sum\limits_{j=0}^{%
	\infty }2^{js\beta }\big\|\mathcal{F}^{-1}\varphi _{j}\ast f\big\|_{L^{p}(%
	\mathbb{R}^{n},w)}^{\beta }\Big)^{1/\beta }
\end{equation*}%
is finite. In the limiting case $\beta =\infty $ the usual modification is
required.

Let $w\in \mathcal{A}_{\infty }$, $s\in \mathbb{R}$, $0<\beta \leq \infty $
and $0<p<\infty $. We define weighted Triebel-Lizorkin spaces $F_{p,\beta
}^{s}(\mathbb{R}^{n},w)$ to be the set of all distributions $f\in \mathcal{S}%
^{\prime }(\mathbb{R}^{n})$ such that%
\begin{equation*}
\big\|f\big\|_{F_{p,\beta }^{s}(\mathbb{R}^{n},w)}=\Big\|\Big(%
\sum\limits_{j=0}^{\infty }2^{js\beta }\left\vert \mathcal{F}^{-1}\varphi
_{j}\ast f\right\vert ^{\beta }\Big)^{1/\beta }\Big\|_{L^{p}(\mathbb{R}%
	^{n},w)}
\end{equation*}%
is finite. In the limiting case $\beta =\infty $ the usual modification is
required.

The spaces $B_{p,\beta }^{s}(\mathbb{R}^{n},w)=B_{p,\beta }^{s}(w)$ and $%
F_{p,\beta }^{s}(\mathbb{R}^{n},w)=F_{p,\beta }^{s}(w)$ are independent of
the particular choice of the smooth dyadic resolution of unity $\{\varphi
_{j}\}_{j\in \mathbb{N}_{0}}$ appearing in their definitions. They are
quasi-Banach spaces (Banach spaces for $p,\beta\geq 1$). Moreover, for $w\equiv
1\in \mathcal{A}_{\infty }$ we obtain the usual (unweighted)
Besov and Triebel-Lizorkin spaces. We refer, in particular, to the papers \cite{Bui82,Bui84,IzSa12} for a comprehensive treatment of the
weighted spaces. Let $w_{\gamma }$ be a power weight, i.e., $w_{\gamma
}(x)=|x|^{\gamma }$ with $\gamma >-n$. Then we have%
\begin{equation*}
B_{p,\beta }^{s}(w_{\gamma })=\dot{K}_{p}^{\frac{\gamma }{p},p}B_{\beta
}^{s}\quad \text{and}\quad F_{p,\beta }^{s}(w_{\gamma })=\dot{K}_{p}^{\frac{%
		\gamma }{p},p}F_{\beta }^{s},
\end{equation*}%
in the sense of equivalent quasi-norms.

\begin{definition}
	{\rm	$\mathrm{(i)}$ Let $1<q<\infty ,0<p<\infty ,-\frac{n}{q}<\alpha <n(1-\frac{1%
		}{q})$ and $s\in \mathbb{R}$. Then the Herz-type Bessel potential space $%
		\dot{k}_{q,s}^{\alpha ,p}$ is the\textit{\ }collection of all $f\in \mathcal{%
			S}^{\prime }(\mathbb{R}^{n})$\ such that%
		\begin{equation*}
		\big\|f\big\|_{\dot{k}_{q,s}^{\alpha ,p}}=\big\|(1+|\xi |^{2})^{\frac{s}{2}%
		}\ast f\big\|_{\dot{K}_{q}^{\alpha ,p}}<\infty .
		\end{equation*}%
		$\mathrm{(ii)}$ Let $1<q<\infty ,0<p<\infty ,-\frac{n}{q}<\alpha <n(1-\frac{1%
		}{q})$ and $m\in \mathbb{N}$. The homogeneous Herz-type Sobolev space $\dot{W%
		}_{q,m}^{\alpha ,p}$ is the collection of all $f\in \mathcal{S}^{\prime }(%
		\mathbb{R}^{n})$\ such that%
		\begin{equation*}
		\big\|f\big\|_{\dot{W}_{q,m}^{\alpha ,p}}=\sum\limits_{|\beta |\leq m}\Big\|%
		\frac{\partial ^{\beta }f}{\partial ^{\beta }x}\Big\|_{\dot{K}_{q}^{\alpha
				,p}}<\infty ,
		\end{equation*}%
		where the derivatives must be understood in the sense of distribution.}
\end{definition}

In the following, we will present the connection between the Herz-type
Triebel-Lizorkin spaces and the Herz-type Bessel potential spaces; see \cite%
{LuYang97, XuYang03}. Let $1<q<\infty ,1<p<\infty $ and $-\frac{n}{%
	q}<\alpha <n(1-\frac{1}{q})$.$\ $If $s\in \mathbb{R}$, then%
\begin{equation}
\dot{K}_{q}^{\alpha ,p}F_{2}^{s}=\dot{k}_{q,s}^{\alpha ,p}
\label{coincidence1}
\end{equation}%
with equivalent norms. If $s=m\in \mathbb{N}$, then

\begin{equation}
\dot{K}_{q}^{\alpha ,p}F_{2}^{m}=\dot{W}_{q,m}^{\alpha ,p}
\label{coincidence2}
\end{equation}%
with equivalent norms. In particular%
\begin{equation*}
\dot{K}_{p}^{0,p}F_{2}^{m}=W_{m}^{p}\quad \text{(Sobolev spaces)}
\end{equation*}%
and%
\begin{equation}
\dot{K}_{q}^{\alpha ,p}F_{2}^{0}=\dot{K}_{q}^{\alpha ,p}
\label{coincidence3}
\end{equation}%
with equivalent norms. Let $0<\theta <1$,%
\begin{equation*}
\alpha =(1-\theta )\alpha _{0}+\theta \alpha _{1},\quad \frac{1}{p}=\frac{%
	1-\theta }{p_{0}}+\frac{\theta }{p_{1}},\quad \frac{1}{q}=\frac{1-\theta }{%
	q_{0}}+\frac{\theta }{q_{1}},\quad \frac{1}{\beta }=\frac{1-\theta }{\beta
	_{0}}+\frac{\theta }{\beta _{1}}
\end{equation*}%
and%
\begin{equation*}
s=(1-\theta )s_{0}+\theta s_{1}.
\end{equation*}%
For simplicity, in what follows, we use $\dot{K}_{q}^{\alpha ,p}A_{\beta
}^{s}$ to denote either $\dot{K}_{q}^{\alpha ,p}B_{\beta }^{s}$ or $\dot{K}%
_{q}^{\alpha ,p}F_{\beta }^{s}$. As an immediate consequence of H\"{o}lder's
inequality we have the so-called interpolation inequalities:%
\begin{equation}
\big\|f\big\|_{\dot{K}_{q}^{\alpha ,p}A_{\beta }^{s}}\leq \big\|f\big\|_{%
	\dot{K}_{q_{0}}^{\alpha _{0},p_{0}}A_{\beta _{0}}^{s_{0}}}^{1-\theta }\big\|f%
\big\|_{\dot{K}_{q_{1}}^{\alpha _{1},p_{1}}A_{\beta _{1}}^{s_{1}}}^{\theta }
\label{Interpolation}
\end{equation}%
holds for all $f\in \dot{K}_{q_{0}}^{\alpha _{0},p_{0}}A_{\beta
	_{0}}^{s_{0}}\cap \dot{K}_{q_{1}}^{\alpha _{1},p_{1}}A_{\beta _{1}}^{s_{1}}$.

We collect some embeddings on these functions spaces as obtained in \cite%
{Drihem1.13}-\cite{Drihem2.13}. First we have elementary embeddings
within these spaces. Let $s\in \mathbb{R},0<p,q<\infty ,0<\beta \leq \infty $%
\ and $\alpha >-\frac{n}{q}$. Then%
\begin{equation}
\dot{K}_{q}^{\alpha ,p}B_{\min \left( \beta ,p,q\right) }^{s}\hookrightarrow 
\dot{K}_{q}^{\alpha ,p}F_{\beta }^{s}\hookrightarrow \dot{K}_{q}^{\alpha
	,p}B_{\max \left( \beta ,p,q\right) }^{s}.  \label{aux7}
\end{equation}

\begin{theorem}
	\label{embeddings3}\textit{Let }$\alpha _{1},\alpha _{2},s_{1},s_{2}\in 
	\mathbb{R},0<s,p,q,r,\beta \leq \infty ,\alpha _{1}>-\frac{n}{s}\ $\textit{%
		and }$\alpha _{2}>-\frac{n}{q}$. \textit{We suppose that }%
	\begin{equation*}
	s_{1}-\frac{n}{s}-\alpha _{1}=s_{2}-\frac{n}{q}-\alpha _{2}.
	\end{equation*}%
	\textit{Let }$0<q\leq s\leq \infty $ and $\alpha _{2}\geq \alpha _{1}$\ or\ $%
	0<s\leq q\leq \infty $ and%
	\begin{equation}
	\alpha _{2}+\frac{n}{q}\geq \alpha _{1}+\frac{n}{s}.  \label{cond-embedding}
	\end{equation}
	
	\noindent $\mathrm{(i)}$\textrm{\ }We have the embedding%
	\begin{equation*}
	\dot{K}_{q}^{\alpha _{2},\theta }B_{\beta }^{s_{2}}\hookrightarrow \dot{K}%
	_{s}^{\alpha _{1},r}B_{\beta }^{s_{1}},
	\end{equation*}%
	where%
	\begin{equation*}
	\theta =\left\{ 
	\begin{array}{ccc}
	r, & \text{if} & \alpha _{2}+\frac{n}{q}=\alpha _{1}+\frac{n}{s},s\leq q%
	\text{ or }\alpha _{2}=\alpha _{1},q\leq s, \\ 
	p, & \text{if} & \alpha _{2}+\frac{n}{q}>\alpha _{1}+\frac{n}{s},s\leq q%
	\text{ or }\alpha _{2}>\alpha _{1},q\leq s.%
	\end{array}%
	\right.
	\end{equation*}
	
	\noindent $\mathrm{(ii)}$ \textit{Let }$0<q,s<\infty $. The embedding%
	\begin{equation*}
	\dot{K}_{q}^{\alpha _{2},r}F_{\theta }^{s_{2}}\hookrightarrow \dot{K}%
	_{s}^{\alpha _{1},p}F_{\beta }^{s_{1}}
	\end{equation*}%
	holds if $0<r\leq p<\infty $, where%
	\begin{equation*}
	\theta =\left\{ 
	\begin{array}{ccc}
	\beta , & \text{if} & 0<s\leq q<\infty \text{ and }\alpha _{2}+\frac{n}{q}%
	=\alpha _{1}+\frac{n}{s}; \\ 
	\infty , &  & \text{otherwise.}%
	\end{array}%
	\right.
	\end{equation*}
\end{theorem}

We now present an immediate consequence of the Sobolev embeddings, which
called Hardy-Sobolev inequalities.

\begin{corollary}
	\textit{Let }$1<q\leq s<\infty $, $1<q<n $ and $\alpha =\frac{n}{q}-\frac{n}{s}-1$.
	There is a constant $c>0$ such that for all $f\in \dot{W}_{q}^{1}$%
	\begin{equation*}
	\int_{\mathbb{R}^{n}}\Big(\frac{|f(x)|}{|x|^{-\alpha }}\Big)^{s}dx\leq c\Big(%
	\sum\limits_{|\beta |=1}\Big\|\frac{\partial ^{\beta }f}{\partial ^{\beta }x}%
	\Big\|_{\dot{K}_{q}^{0,s}}\Big)^{s}\leq c\Big(\sum\limits_{|\beta |=1}\Big\|%
	\frac{\partial ^{\beta }f}{\partial ^{\beta }x}\Big\|_{q}\Big)^{s}.
	\end{equation*}
\end{corollary}

Now we recall the Franke embedding, see \cite{drihem2016jawerth}.

\begin{theorem}
	\label{Franke-emb}\textit{Let }$\alpha _{1},\alpha _{2},s_{1},s_{2}\in 
	\mathbb{R}$, $0<s,p,q<\infty ,0<\theta \leq \infty ,\alpha _{1}>-\frac{n}{s}%
	\ $\textit{and }$\alpha _{2}>-\frac{n}{q}$. \textit{We suppose that }%
	\begin{equation*}
	s_{1}-\frac{n}{s}-\alpha _{1}=s_{2}-\frac{n}{q}-\alpha _{2}.
	\end{equation*}
	
	\noindent \textit{Let }%
	\begin{equation*}
	0<q<s<\infty \quad \text{and}\quad \alpha _{2}\geq \alpha _{1},
	\end{equation*}%
	or 
	\begin{equation*}
	0<s\leq q<\infty \quad \text{and}\quad \alpha _{2}+\frac{n}{q}>\alpha _{1}+%
	\frac{n}{s}.
	\end{equation*}%
	Then%
	\begin{equation*}
	\dot{K}_{q}^{\alpha _{2},p}B_{p}^{s_{2}}\hookrightarrow \dot{K}_{s}^{\alpha
		_{1},p}F_{\theta }^{s_{1}}.
	\end{equation*}
\end{theorem}

\begin{corollary}
	\textit{Let }$1<q\leq s<\infty $ with $1<q<n $. \textit{Let }$\alpha =\frac{n}{q}-\frac{n}{s}-1$. There is a constant $c>0$ such that for all $f\in B_{q,s}^{1}$%
	\begin{equation*}
	\int_{\mathbb{R}^{n}}\Big(\frac{|f(x)|}{|x|^{-\alpha }}\Big)^{s}dx\leq c%
	\big\|f\big\|_{\dot{K}_{q}^{0,s}B_{s}^{1}}^{s}\leq c\big\|f\big\|%
	_{B_{q,s}^{1}}^{s}.
	\end{equation*}
\end{corollary}

\begin{remark}
	{\rm	We would like to mention that in Theorem \ref{embeddings3} and Theorem \ref%
		{Franke-emb} the assumptions $s_{1}-\tfrac{n}{s}-\alpha _{1}\leq s_{2}-%
		\tfrac{n}{q}-\alpha _{2}$, \eqref{cond-embedding} and $0<r\leq p<\infty $
		are necessary, see \cite{Drihem1.13,Drihem2.13, drihem2016jawerth}.}
\end{remark}

Let  $\{\varphi _{j}\}_{j\in \mathbb{N}_{0}}$ be a resolution of unity. For any $a>0$, $f\in \mathcal{S}^{\prime }(\mathbb{R}^{n})$ and $x\in 
\mathbb{R}^{n}$, we denote, Peetre maximal function, 
\begin{equation*}
(\mathcal{F}^{-1}\varphi _{j})^{\ast ,a}f(x)=\sup_{y\in \mathbb{R}^{n}}\,\frac{\left\vert
	\mathcal{F}^{-1}\varphi _{j}\ast f(y)\right\vert }{(1+2^{j}\left\vert x-y\right\vert )^{a}}%
,\quad j\in \mathbb{N}_{0}.
\end{equation*}%
We now present a fundamental characterization of the above spaces, which
plays an essential role in this paper, see \cite[Theorem 1]{Xu05}.

\begin{theorem}
	\label{Peetre maximal function}Let $s\in \mathbb{R},0<p,q<\infty ,0<\beta
	\leq \infty $ and $\alpha >-\frac{n}{q}$. Let  $a>\frac{n}{\min \big(q,\frac{n}{\alpha +\frac{n}{q}}\big)%
	}$.
	Then 
	
	\begin{equation*}
	\big\Vert f\big\Vert_{\dot{K}_{q}^{\alpha ,p}B_{\beta }^{s}}^{\star }=%
	\Big(\sum_{j=0}^{\infty }2^{js\beta }\big\Vert(\mathcal{F}^{-1}\varphi _{j})^{\ast
		,a}f\big\Vert_{\dot{K}_{q}^{\alpha ,p}}^{\beta }\Big)^{1/\beta },
	\end{equation*}
	is an equivalent quasi-norm in $\dot{K}_{q}^{\alpha ,p}B_{\beta }^{s}.$ Let $a>\frac{n}{\min \big(\min (q,\beta),\frac{n}{\alpha +\frac{n}{q}}\big)}$.
	Then 
	\begin{equation*}
	\big\Vert f\big\Vert_{\dot{K}_{q}^{\alpha ,p}F_{\beta }^{s}}^{\star }=%
	\Big\Vert\Big(\sum_{j=0}^{\infty }2^{js\beta }(\mathcal{F}^{-1}\varphi _{j})^{\ast
		,a}f)^{\beta }\Big)^{1/\beta }\Big\Vert_{\dot{K}_{q}^{\alpha ,p}},
	\end{equation*}
	is an\ equivalent quasi-norm in $\dot{K}_{q}^{\alpha ,p}F_{\beta }^{s}.$
\end{theorem}

Let $0<p,q\leq \infty $. For later use we introduce the following
abbreviations:%
\begin{equation*}
\sigma _{q}=n\max (\frac{1}{q}-1,0)\quad \text{and}\quad \sigma _{p,q}=n\max
(\frac{1}{p}-1,\frac{1}{q}-1,0).
\end{equation*}%
In the next we shall interpret $L_{\mathrm{loc}}^{1}$ as the set of regular
distributions.

\begin{theorem}
	\label{regular-distribution1}\textit{Let }$0<p,q,\beta \leq \infty ,\alpha >-%
	\frac{n}{q}$, $\alpha _{0}=n-\frac{n}{q}$ and $s>\max (\sigma _{q},\alpha
	-\alpha _{0})$. Then%
	\begin{equation*}
	\dot{K}_{q}^{\alpha ,p}A_{\beta }^{s}\hookrightarrow L_{\mathrm{loc}}^{1},
	\end{equation*}%
	where $0<p,q<\infty $ in the case of Herz-type Triebel-Lizorkin spaces.
\end{theorem}

\begin{proof} Let $\{\varphi _{j}\}_{j\in \mathbb{N}_{0}}$ be a smooth dyadic resolution of unity. We set
	\begin{equation*}
	\varrho _{k}=\sum\limits_{j=0}^{k}\mathcal{F}^{-1}\varphi _{j}\ast f,\quad
	k\in \mathbb{N}_{0}.
	\end{equation*} For technical reasons, we split the proof into two steps.
	
	\textit{Step 1.} We consider the case $1\leq q\leq\infty $. In order to prove we
	additionally do it into the four Substeps 1.1, 1.2, 1.3 and 1.4.
	
	\textit{Substep 1.1.} $-\frac{n}{q}<\alpha <\alpha _{0}$.  
	Since $s > 0$  and $\dot{K}_{q}^{\alpha
		,p}\hookrightarrow \dot{K}_{q}^{\alpha
		,\max (1,p)}$, we have 
	\begin{equation*}
	\sum_{j=0}^{\infty }\big\|\mathcal{F}^{-1}\varphi _{j}\ast f\big\|_{\dot{K}%
		_{q}^{\alpha ,\max (1,p)}}\lesssim \big\|f\big\|_{\dot{K}_{q}^{\alpha
			,p}A_{\beta }^{s}}.
	\end{equation*}%
	Then, the sequence $\{\varrho _{k}\}_{k\in \mathbb{N}_{0}}$ converges to $g \in  \dot{K}_{q}^{\alpha ,\max (1,p)}$. Let $\varphi \in \mathcal{S}(\mathbb{R}^{n})$. Write 
	\begin{equation*}
	\langle f-g,\varphi \rangle =\langle f-\varrho _{N},\varphi \rangle +\langle
	g-\varrho _{N},\varphi \rangle ,\quad N\in \mathbb{N}_{0}.
	\end{equation*}%
	Here $\langle \cdot ,\cdot \rangle $ denotes the duality bracket between $%
	\mathcal{S}^{\prime }(\mathbb{R}^{n})$ and $\mathcal{S}(\mathbb{R}^{n})$.
	Clearly, the first term tends to zero as $N\rightarrow \infty $, while by H\"{o}lder's\ inequality there exists a constant $C>0$ independent of $N$
	such that 
	\begin{equation*}
	|\langle g-\varrho _{N},\varphi \rangle |\leq C\big\|g-\varrho _{N}\big\|_{%
		\dot{K}_{q}^{\alpha ,\max (1,p)}},
	\end{equation*}%
	which tends to zero as $N\rightarrow \infty $. 
	From this and $\dot{K}_{q}^{\alpha ,\max (1,p)}\hookrightarrow L_{\mathrm{loc}}^{1}$, because of $\alpha <\alpha _{0}$, we deduce the desired result. In addition, we have
	\begin{equation*}
	\dot{K}_{q}^{\alpha ,p}A_{\beta }^{s}\hookrightarrow \dot{K}_{q}^{\alpha
		,\max (1,p)}.
	\end{equation*}%
	
	\textit{Substep 1.2.} $\alpha \geq \alpha _{0}$ and $1<q\leq\infty $.
	Let $1<q_{1}<\infty $ be such that%
	\begin{equation*}
	s>\alpha +\frac{n}{q}-\frac{n}{q_{1}}.
	\end{equation*}%
	We distinguish two cases:
	
	$\bullet $ $q_{1}=q$. By Theorem \ref{embeddings3}/(i), we obtain%
	\begin{equation*}
	\dot{K}_{q}^{\alpha ,p}B_{\beta }^{s}\hookrightarrow \dot{K}%
	_{q}^{0,q}B_{\beta }^{s-\alpha }=B_{q,\beta }^{s-\alpha }\hookrightarrow L_{%
		\mathrm{loc}}^{1}.
	\end{equation*}%
	where the last embedding follows by the fact that%
	\begin{equation}
	B_{q,\beta }^{s-\alpha }\hookrightarrow L^{q},  \label{Substep1.2.1}
	\end{equation}%
	because of $s-\alpha >0$. The Herz-type Triebel-Lizorkin case follows by the
	second embeddings of \eqref{aux7}.
	
	$\bullet $ $1<q_{1}<q\leq\infty $ or $1<q<q_{1}<\infty $. If we
	assume the first possibility then Theorem \ref{embeddings3}/(i) and Substep 1.1 yield%
	\begin{equation*}
	\dot{K}_{q}^{\alpha ,p}B_{\beta }^{s}\hookrightarrow\dot{K}_{q_{1}}^{0,p}B_{\beta }^{s-\alpha-\frac{n}{q}+\frac{n}{q_{1}} } \hookrightarrow
	L_{\mathrm{loc}}^{1},
	\end{equation*}%
	since  $\alpha +\frac{n}{q}>\frac{n}{q_{1}}$. The latter possibility follows again by Theorem \ref{embeddings3}/(i). Indeed, we have
	\begin{equation*}
	\dot{K}_{q}^{\alpha ,p}B_{\beta }^{s }\hookrightarrow\dot{K}_{q}^{\alpha _{0},p}B_{\beta }^{s+\alpha _{0}-\alpha }\hookrightarrow 
	\dot{K}_{q_{1}}^{0,q_{1}}B_{\beta }^{s-\alpha -\frac{n}{q}+\frac{n}{q_{1}}%
	}=B_{q_{1},\beta }^{s-\alpha -\frac{n}{q}+\frac{n}{q_{1}}}\hookrightarrow
	L_{\mathrm{loc}}^{1},
	\end{equation*}%
	where the last embedding follows\ by the fact that%
	\begin{equation}
	B_{q_{1},\beta }^{s-\alpha -\frac{n}{q}+\frac{n}{q_{1}}}\hookrightarrow
	L^{q_{1}}.  \label{Substep1.2.2}
	\end{equation}%
	Therefore from \eqref{aux7} we obtain the desired embeddings.
	
	\textit{Substep 1.3.} $q=1$ and $\alpha >0$. We have%
	\begin{equation*}
	\dot{K}_{1}^{\alpha ,p}B_{\beta }^{s}\hookrightarrow \dot{K}%
	_{1}^{0,1}B_{\beta }^{s-\alpha }=B_{1,\beta }^{s-\alpha }\hookrightarrow
	L^{1},
	\end{equation*}%
	since $s>\alpha $.
	
	\textit{Substep 1.4.} $q=1$ and $\alpha =0$.  Let $\alpha _{3}$ be a real
	number such that $\max (-n,-s)<\alpha _{3}<0$. From Theorem \ref{embeddings3}%
	, we get 
	\begin{equation*}
	\dot{K}_{1}^{0,p}A_{\beta }^{s}\hookrightarrow \dot{K}_{1}^{\alpha
		_{3},p}A_{\beta }^{s+\alpha _{3}}.
	\end{equation*}%
	We have 
	\begin{equation*}
	\sum_{k=0}^{\infty }\big\|\mathcal{F}^{-1}\varphi _{k}\ast f\big\|_{\dot{K}%
		_{1}^{\alpha _{3},\max (1,p)}}\lesssim \big\|f\big\|_{\dot{K}_{1}^{\alpha
			_{3},p}A_{\beta }^{s+\alpha _{3}}}\lesssim \big\|f\big\|_{\dot{K}%
		_{1}^{0,p}A_{\beta }^{s}},
	\end{equation*}%
	since $\alpha _{3}+s>0$. Using the same type of arguments as in Substep 1.1 it is easy to see that
	\begin{equation*}
	\dot{K}_{1}^{\alpha _{3},p}A_{\beta }^{s+\alpha _{3}}\hookrightarrow \dot{K}%
	_{1}^{\alpha _{3},\max (1,p)}\hookrightarrow L_{\mathrm{loc}}^{1}.
	\end{equation*}
	
	\textit{Step 2.} We consider the case $0<q<1$.
	
	\textit{Substep 2.1.} $-\frac{n}{q}<\alpha <0$. By Lemma \ref%
	{Bernstein-Herz-ine1}, we obtain%
	\begin{equation*}
	\sum_{j=0}^{\infty }\big\|\mathcal{F}^{-1}\varphi _{j}\ast f\big\|_{\dot{K}%
		_{1}^{\alpha ,\max (1,p)}}\lesssim \sum_{j=0}^{\infty }2^{j(\frac{n}{q}-n)}%
	\big\|\mathcal{F}^{-1}\varphi _{j}\ast f\big\|_{\dot{K}_{q}^{\alpha
			,p}}\lesssim \big\|f\big\|_{\dot{K}_{q}^{\alpha ,p}A_{\beta }^{s}},
	\end{equation*}%
	since $s>\frac{n}{q}-n$.  The desired embedding follows by the fact that $\dot{K}_{1}^{\alpha ,\max
		(1,p)}\hookrightarrow L_{\mathrm{loc}}^{1}$ and the arguments in Substep 1.1.  In addition%
	\begin{equation}
	\dot{K}_{q}^{\alpha ,p}A_{\beta }^{s}\hookrightarrow \dot{K}_{1}^{\alpha
		,\max (1,p)}.  \label{q-less1}
	\end{equation}%

	\textit{Substep 2.2.} $\alpha \geq 0$. Let $\alpha _{4}$ be a real
	number such that $\max (-n,-s+\frac{n}{q}-n+\alpha)<\alpha _{4}<0$. From Theorem \ref{embeddings3}%
	, we get 
	\begin{equation*}
	\dot{K}_{q}^{\alpha,p}A_{\beta }^{s}\hookrightarrow\dot{K}_{1}^{0,p}A_{\beta }^{s-\frac{n}{q}+n-\alpha}\hookrightarrow \dot{K}_{1}^{\alpha
		_{4},p}A_{\beta }^{s-\frac{n}{q}+n-\alpha+\alpha _{4}}\hookrightarrow\dot{K}_{1}^{\alpha
		_{4},\max(1,p)}A_{\beta }^{s-\frac{n}{q}+n-\alpha+\alpha _{4}} \label{alpha=0}.
	\end{equation*}
	As in Substep 1.4, we easily obtain that%
	\begin{equation*}
	\dot{K}_{q}^{\alpha ,p}A_{\beta }^{s}\hookrightarrow \hookrightarrow L_{\mathrm{loc}%
	}^{1}.  \label{alpha-positive}
	\end{equation*}
	Therefore, under the hypothesis of this theorem, every $f\in \dot{K}%
	_{q}^{\alpha ,p}A_{\beta }^{s}$ is a regular distribution. This finishes the
	proof.
\end{proof}

Let $f$ be an arbitrary function on $\mathbb{R}^{n}$ and $x,h\in \mathbb{R}%
^{n}$. Then%
\begin{equation*}
\Delta _{h}f(x)=f(x+h)-f(x),\quad \Delta _{h}^{M+1}f(x)=\Delta _{h}(\Delta
_{h}^{M}f)(x),\quad M\in \mathbb{N}.
\end{equation*}%
These are the well-known differences of functions which play an important
role in the theory of function spaces. Using mathematical induction one can
show the explicit formula%
\begin{equation*}
\Delta _{h}^{M}f(x)=\sum_{j=0}^{M}\left( -1\right) ^{j}C_{j}^{M}f(x+(M-j)h),\quad x \in \mathbb{R}^{n},
\end{equation*}  %
where $C_{j}^{M}$ are the binomial coefficients. By ball means of
differences we mean the quantity%
\begin{equation*}
d_{t}^{M}f(x)=t^{-n}\int_{|h|\leq t}\left\vert \Delta
_{h}^{M}f(x)\right\vert dh=\int_{B}\left\vert \Delta
_{th}^{M}f(x)\right\vert dh,\quad x \in \mathbb{R}^{n}.
\end{equation*}%
Here $B=\{y\in \mathbb{R}^{n}:|h|\leq 1\}$ is the unit ball of $\mathbb{R}%
^{n}$ and $t>0$ is a real number. We set%
\begin{equation*}
\big\|f\big\|_{\dot{K}_{q}^{\alpha ,p}B_{\beta }^{s}}^{\ast }=\big\|f\big\|_{%
	\dot{K}_{q}^{\alpha ,p}}+\Big(\int_{0}^{\infty }t^{-s\beta }\big\|d_{t}^{M}f%
\big\|_{\dot{K}_{q}^{\alpha ,p}}^{\beta }\frac{dt}{t}\Big)^{1/\beta }
\label{norm1}
\end{equation*}%
and%
\begin{equation*}
\big\|f\big\|_{\dot{K}_{q}^{\alpha ,p}F_{\beta }^{s}}^{\ast }=\big\|f\big\|_{%
	\dot{K}_{q}^{\alpha ,p}}+\Big\|\Big(\int_{0}^{\infty }t^{-s\beta
}(d_{t}^{M}f)^{\beta }\frac{dt}{t}\Big)^{1/\beta }\Big\|_{\dot{K}%
	_{q}^{\alpha ,p}}.  \label{norm2}
\end{equation*}%
The following theorem play a central role in our paper.

\begin{theorem}
	\label{means-diff-cha}\textit{Let }$0<p,q,\beta \leq \infty ,\alpha >-\frac{n%
	}{q}$, $\alpha _{0}=n-\frac{n}{q}\ $and $M\in \mathbb{N}\backslash \{0\}.%
	\newline
	\mathrm{(i)}$ Assume that%
	\begin{equation*}
	\max (\sigma _{q},\alpha -\alpha _{0})<s<M.
	\end{equation*}%
	Then $\left\Vert \cdot \right\Vert _{\dot{K}_{q}^{\alpha ,p}B_{\beta
		}^{s}}^{\ast }$ is an equivalent quasi-norm on $\dot{K}_{q}^{\alpha
		,p}B_{\beta }^{s}$.$\newline
	\mathrm{(ii)}$ Let\ $0<p<\infty \ $and $0<q<\infty $. Assume that%
	\begin{equation*}
	\max (\sigma _{q,\beta },\alpha -\alpha _{0})<s<M.
	\end{equation*}%
	Then $\left\Vert \cdot \right\Vert _{\dot{K}_{q}^{\alpha ,p}F_{\beta
		}^{s}}^{\ast }$ is an equivalent quasi-norm on $\dot{K}_{q}^{\alpha
		,p}F_{\beta }^{s}$.
\end{theorem}

\begin{proof}
	For ease of presentation, we split the proof into three steps.
	
	\textit{Step 1}. We will prove that
	\begin{equation*}
	\big\|f\big\|_{\dot{K}_{q}^{\alpha ,p}}\lesssim \big\|f\big\|_{\dot{K}%
		_{q}^{\alpha ,p}A_{\beta }^{s}}
	\end{equation*}%
	for all $f\in \dot{K}_{q}^{\alpha ,p}A_{\beta }^{s}$. We employ the same
	notations as in Theorem \ref{regular-distribution1}. Recall that 
	\begin{equation*}
	\varrho _{k}=\sum\limits_{j=0}^{k}\mathcal{F}^{-1}\varphi _{j}\ast f,\quad
	k\in \mathbb{N}_{0}.
	\end{equation*}%
	Obviously $\{\varrho _{k}\}_{k\in \mathbb{N}_{0}}$ converges to $f$ in $\mathcal{S}^{\prime }(%
	\mathbb{R}^{n})$ and $\{\varrho _{k}\}_{k\in \mathbb{N}_{0}}\subset \dot{K}_{q}^{\alpha ,p}$ for
	any $0<p,q\leq \infty \ $and any $\alpha >-\frac{n}{q}$. Furthermore, $%
	\{\varrho _{k}\}_{k\in \mathbb{N}_{0}}$ is a Cauchy sequences in $\dot{K}_{q}^{\alpha ,p}$ and
	hence it converges to a function  $g\in \dot{K}_{q}^{\alpha ,p}$, and%
	\begin{equation*}
	\big\|g\big\|_{\dot{K}_{q}^{\alpha ,p}}\lesssim \big\|f\big\|_{\dot{K}%
		_{q}^{\alpha ,p}A_{\beta }^{s}}.
	\end{equation*}%
	Let us prove that $g=f$ a.e. We will do this into four cases.
	
	\noindent \textit{Case 1.} $-\frac{n}{q}<\alpha <\alpha _{0}$ and $1\leq
	q\leq\infty $. Let $\varphi \in \mathcal{D}(\mathbb{R}^{n})$. We write 
	\begin{equation*}
	\langle f-g,\varphi \rangle =\langle f-\varrho _{N},\varphi \rangle +\langle
	g-\varrho _{N},\varphi \rangle ,\quad N\in \mathbb{N}_{0}.
	\end{equation*}%
	Here $\langle \cdot ,\cdot \rangle $ denotes the duality bracket between $%
	\mathcal{S}^{\prime }(\mathbb{R}^{n})$ and $\mathcal{S}(\mathbb{R}^{n})$.
	Clearly, the first term tends to zero as $N\rightarrow \infty $, while by H%
	\"{o}lder's\ inequality there exists a constant $C>0$ independent of $N$
	such that 
	\begin{equation*}
	|\langle g-\varrho _{N},\varphi \rangle |\leq C\big\|g-\varrho _{N}\big\|_{%
		\dot{K}_{q}^{\alpha ,\max (1,p)}},
	\end{equation*}%
	which tends to zero as $N\rightarrow \infty $. Then, with the help of Substep 1.1 of the proof of Theorem \ref{regular-distribution1}, we have  $g=f$ almost
	everywhere.
	
	\textit{Case 2.} $\alpha \geq \alpha _{0}$ and $1<q\leq\infty $. Let $%
	1<q_{1}<\infty $ be as in Theorem \ref{regular-distribution1}. From %
	\eqref{Substep1.2.1} and \eqref{Substep1.2.2}, we derive in this case, that
	every $f\in \dot{K}_{q}^{\alpha ,p}A_{\beta }^{s}$ is a regular
	distribution, $\{\varrho _{k}\}_{k\in \mathbb{N}_{0}}$ converges to $f$ in $L^{q_{1}}$ and 
	\begin{equation*}
	\big\|f\big\|_{q_{1}}\lesssim \big\|f\big\|_{\dot{K}_{q}^{\alpha ,p}A_{\beta
		}^{s}}.
	\end{equation*}%
	Indeed, from the embeddings \eqref{Substep1.2.2} and since $f\in
	B_{q_{1},\beta }^{\frac{n}{q_{1}}-\alpha -\frac{n}{q}+s}$, it follows that $\{\varrho _{k}\}_{k\in \mathbb{N}_{0}}$ converges to a function $h$ $ \in$ $L^{q_{1}}$. Similarly as in Case 1,
	we conclude that  $f=h$ a.e. It remains to prove that $g=f$ a.e. We have 
	\begin{equation*}
	\big\|f-g\big\|_{\dot{K}_{q}^{\alpha ,p}}^{\sigma }\leq \big\|f-\varrho _{k}%
	\big\|_{\dot{K}_{q}^{\alpha ,p}}^{\sigma }+\big\|g-\varrho _{k}\big\|_{\dot{K%
		}_{q}^{\alpha ,p}}^{\sigma },\quad k\in \mathbb{N}_{0}
	\end{equation*}%
	and%
	\begin{equation*}
	\big\|f-\varrho _{k}\big\|_{\dot{K}_{q}^{\alpha ,p}}^{\sigma } \leq
	\sum\limits_{j=k+1}^{\infty }\big\|\mathcal{F}^{-1}\varphi _{j}\ast f\big\|%
	_{\dot{K}_{q}^{\alpha ,p}}^{\sigma } \leq \big\|f\big\|_{\dot{K}_{q}^{\alpha ,p}A_{\beta }^{s}}^{\sigma
	}\sum\limits_{j=k+1}^{\infty }2^{-js\sigma },
	\end{equation*}%
	where $\sigma =\min (1,p,q)$. Letting $k$ tends to infinity, we get $g=f$ a.e.  For the latter case $1<q_{1}<q\leq\infty $, we have
	\begin{equation*}
	\dot{K}_{q}^{\alpha,p}A_{\beta }^{s}\hookrightarrow \dot{K}_{q_{1}}^{0,\max(1,p)}A_{\beta }^{s-\alpha-\frac{n}{q}+\frac{n}{q_{1}}}.
	\end{equation*}
	As in Case 1, $\{\varrho _{k}\}_{k\in \mathbb{N}_{0}}$ converges to a function $h$ $ \in$ $\dot{K}_{q_{1}}^{0,\max(1,p)}$. Then again, similarly to the arguments in Case 1 it is easy to check that  $f=h$ a.e. Therefore, we can conclude that $g=f$ a.e. 
	
	\textit{Case 3.} $q=1$  and $ \alpha\geq0.$
	
	\textit{Subcase 3.1. } $q=1$  and $\alpha >0$. We have%
	\begin{equation*}
	\dot{K}_{1}^{\alpha ,p}B_{\beta }^{s}\hookrightarrow L^{1},
	\end{equation*}%
	since $s>\alpha $, see Theorem \ref{regular-distribution1}, Substep 1.3. Now one can continue as in Case 2.
	
	\textit{Subcase 3.2. }  $q=1$  and $\alpha =0$.  Let $\alpha _{3}$ be a real
	number such that $\max (-n,-s)<\alpha _{3}<0$. From Theorem \ref{embeddings3}%
	, we get 
	\begin{equation*}
	\dot{K}_{1}^{0,p}A_{\beta }^{s}\hookrightarrow \dot{K}_{1}^{\alpha
		_{3},p}A_{\beta }^{s+\alpha _{3}}.
	\end{equation*}%
	We have 
	\begin{equation*}
	\sum_{k=0}^{\infty }\big\|\mathcal{F}^{-1}\varphi _{k}\ast f\big\|_{\dot{K}%
		_{1}^{\alpha _{3},\max (1,p)}}\lesssim \big\|f\big\|_{\dot{K}_{1}^{\alpha
			_{3},p}A_{\beta }^{s+\alpha _{3}}}\lesssim \big\|f\big\|_{\dot{K}%
		_{1}^{0,p}A_{\beta }^{s}},
	\end{equation*}%
	since $\alpha _{3}+s>0$. Hence the sequence $\{\varrho _{k}\}_{k\in \mathbb{N}_{0}}$ converges to $%
	f$ in $\dot{K}_{1}^{\alpha _{3},\max (1,p)}$, see Case 1. As in Case 2, we
	obtain $g=f$ a.e.
	
	\textit{Case 4.} $0<q<1$.
	
	\textit{Subcase 4.1. }$-\frac{n}{q}<\alpha <0$. From the embedding %
	\eqref{q-less1} and the fact that $s>\frac{n}{q}-n$, the sequence $\{\varrho
	_{k}\}_{k\in \mathbb{N}_{0}}$ converge to $f$ in $\dot{K}_{1}^{\alpha ,\max (1,p)}$. As above we
	prove that $g=f$ a.e.
	
	\textit{Subcase 4.2. }$\alpha \geq 0$. Recall that 
	\begin{equation*}
	\dot{K}_{q}^{\alpha,p}A_{\beta }^{s}\hookrightarrow \dot{K}_{1}^{\alpha
		_{4},\max(1,p)}A_{\beta }^{s-\frac{n}{q}+n-\alpha+\alpha _{4}} \label{alpha=0},
	\end{equation*}
	see Substep 2.2 of the proof of Theorem \ref{regular-distribution1}. As in Subcase 3.2 the sequence $\{\varrho _{k}\}_{k\in \mathbb{N}_{0}}$ converges to $f$ in $\dot{K}_{1}^{\alpha _{4},\max(1,p) }$. The same arguments above one can
	conclude that: $g=f$ a.e..
	
	\textit{Step 2.} In this step we prove that%
	\begin{equation*}
	\big\|f\big\|_{\dot{K}_{q}^{\alpha ,p}F_{\beta }^{s}}^{\ast \ast }=\Big\|%
	\Big(\int_{0}^{\infty }t^{-s\beta }(d_{t}^{M}f)^{\beta }\frac{dt}{t}\Big)^{%
		1/\beta }\Big\|_{\dot{K}_{q}^{\alpha ,p}}\lesssim \big\|f\big\|_{%
		\dot{K}_{q}^{\alpha ,p}F_{\beta }^{s}},\quad f\in \dot{K}_{q}^{\alpha
		,p}F_{\beta }^{s}.
	\end{equation*}%
	Thus, we need to prove that%
	\begin{equation*}
	\Big\|\Big(\sum_{k=-\infty }^{\infty }2^{sk\beta }|d_{2^{-k}}^{M}f|^{\beta }%
	\Big)^{1/\beta }\Big\|_{\dot{K}_{q}^{\alpha,p }}
	\end{equation*}%
	does not exceed $c\big\|f\big\|_{\dot{K}_{q}^{\alpha ,p}F_{\beta }^{s}}$. In
	order to prove we additionally do it into the two  Substeps 2.1 and
	2.2. The estimate for the space $\dot{K}_{q}^{\alpha ,p}B_{\beta }^{s}$ is
	similar.
	
	\textit{Substep 2.1.} We will estimate%
	\begin{equation*}
	\Big\|\Big(\sum_{k=0}^{\infty }2^{sk\beta }|d_{2^{-k}}^{M}f|^{\beta }\Big)%
	^{1/\beta }\Big\|_{\dot{K}_{q}^{\alpha,p }}.
	\end{equation*}%
	Let $\{ \varphi _{j}\}_{j\in \mathbb{N}_{0}}$\ be a smooth dyadic
	resolution of unity. Obviously we need to estimate%
	\begin{equation}
	\Big\{2^{ks}\sum\limits_{j=0}^{k}d_{2^{-k}}^{M}(\mathcal{F}^{-1}\varphi _{j}\ast f)\Big\}%
	_{k\in \mathbb{N}_{0}}  \label{First-term}
	\end{equation}%
	and 
	\begin{equation}
	\Big\{2^{ks}\sum\limits_{j=k+1}^{\infty }d_{2^{-k}}^{M}(\mathcal{F}^{-1}\varphi _{j}\ast f)%
	\Big\}_{k\in \mathbb{N}_{0}}.  \label{second-term}
	\end{equation}%
	Recall that%
	\begin{equation*}
	d_{2^{-k}}^{M}(\mathcal{F}^{-1}\varphi _{j}\ast f)\lesssim 2^{\left( j-k\right) M}(\mathcal{F}^{-1}\varphi_{j})
	^{\ast ,a}f\left( x\right) 
	\end{equation*}%
	if $a>0$, $0\leq j\leq k,k\in \mathbb{N}_{0}$ and $x\in \mathbb{R}^{n}$, see, e.g., 
	\cite{D22}, where the implicit constant is independent of $j,k$ and $x$.
	We choose $a>\frac{n}{\min \big(\min (q,\beta),\frac{n}{\alpha +\frac{n}{q}}\big)}$. Since $s<M$, \eqref{First-term} in $\ell ^{\beta }$-quasi-norm does not exceed%
	\begin{equation}
	\Big(\sum\limits_{j=0}^{\infty }2^{js\beta }((\mathcal{F}^{-1}\varphi _{j})^{\ast ,a}f)^{\beta
	}\Big)^{1/\beta }.  \label{First-term1}
	\end{equation}%
	By Theorem \ref{Peetre maximal function}, the $\dot{K}_{q}^{\alpha ,p}$%
	-quasi-norm of \eqref{First-term1} is bounded by $c\big\|f\big\|_{\dot{K}%
		_{q}^{\alpha ,p}F_{\beta }^{s}}.$
	
	Now, we estimate \eqref{second-term}. We can distinguish two cases as
	follows:
	
	$\bullet $ \textit{Case 1.  }$\min (q,\beta )\leq 1$. If $-\frac{n}{q}%
	<\alpha <n(1-\frac{1}{q})$, then $s>\frac{n}{\min (q,\beta )}-n$. We choose 
	\begin{equation}
	\max \Big(0,1-\frac{s\min (q,\beta )}{n}\Big)<\lambda <\min (q,\beta ),
	\label{lamda5}
	\end{equation}%
	which is possible because of 
	\begin{equation*}
	s>\frac{n}{\min (q,\beta )}-n=\frac{n}{\min (q,\beta )}\Big(1-\min (q,\beta )%
	\Big).
	\end{equation*}%
	Let $\frac{n}{\min (q,\beta )}<a<\frac{s}{1-\lambda }$. Then $s>a(1-\lambda )
	$. Now, assume that $\alpha \geq n(1-\frac{1}{q})$. Therefore 
	\begin{equation*}
	s>\max \Big(\frac{n}{\min (q,\beta )}-n,\frac{n}{q}+\alpha -n\Big).
	\end{equation*}%
	If $\min (q,\beta )\leq \frac{n}{\frac{n}{q}+\alpha }$, then we choose $%
	\lambda $ as in \eqref{lamda5}. If $\min (q,\beta )>\frac{n}{\frac{n}{q}%
		+\alpha }$, then we choose  
	\begin{equation}
	\max \Big(0,1-\frac{s}{\frac{n}{q}+\alpha }\Big)<\lambda <\frac{n}{\frac{n}{q%
		}+\alpha }  \label{lamda6}
	\end{equation}%
	which is possible because of 
	\begin{equation*}
	s>\frac{n}{q}+\alpha -n=\big(\frac{n}{q}+\alpha \big)\big(1-\frac{n}{\frac{n%
		}{q}+\alpha }\big).
	\end{equation*}%
	In that case, we choose $\frac{n}{q}+\alpha <a<\frac{s}{1-\lambda }$. We set 
	\begin{equation*}
	J_{2,k}(f)=2^{ks}\sum\limits_{j=k+1}^{\infty }d_{2^{-k}}^{M}(\mathcal{F}%
	^{-1}\phi _{j}\ast f),\quad k\in \mathbb{N}_{0}.
	\end{equation*}%
	Recalling the definition of $d_{2^{-k}}^{M}(\phi _{j}\ast f)$, we have 
	\begin{align}
	d_{2^{-k}}^{M}(\mathcal{F}^{-1}\phi _{j}\ast f)& =\int_{B}\big|\Delta
	_{2^{-k}h}^{M}(\mathcal{F}^{-1}\phi _{j}\ast f)\big|dh  \notag \\
	& \leq \int_{B}\big|\Delta _{2^{-k}h}^{M}(\mathcal{F}^{-1}\phi _{j}\ast f)%
	\big|^{\lambda }dh\sup_{h\in B}\big|\Delta _{2^{-k}h}^{M}(\mathcal{F}%
	^{-1}\phi _{j}\ast f)\big|^{1-\lambda }.  \label{term2}
	\end{align}%
	Observe that 
	\begin{equation}
	\big|\mathcal{F}^{-1}\phi _{j}\ast f(x+(M-i)2^{-k}h)\big|\leq c2^{\left(
		j-k\right) a}\phi _{j}^{\ast ,a}f\left( x\right) ,\quad |h|\leq 1
	\label{term3}
	\end{equation}%
	and 
	\begin{equation}
	\int_{B}\big|\mathcal{F}^{-1}\phi _{j}\ast f(x+(M-i)2^{-k}h)\big|^{\lambda
	}dh\leq c\mathcal{M}(|\mathcal{F}^{-1}\phi _{j}\ast f|^{\lambda })(x).
	\label{term4}
	\end{equation}%
	if $j>k,i\in \{0,...,M\}$ and $x\in \mathbb{R}^{n}$. Therefore 
	\begin{equation*}
	d_{2^{-k}}^{M}(\mathcal{F}^{-1}\phi _{j}\ast f)\leq c2^{\left( j-k\right)
		a(1-\lambda )}(\phi _{j}^{\ast ,a}f)^{1-\lambda }\mathcal{M}(|\mathcal{F}%
	^{-1}\phi _{j}\ast f|^{\lambda })
	\end{equation*}%
	for any $j>k$, where the positive constant $c$ is independent of $j$ and $k$%
	. Hence 
	\begin{equation*}
	J_{2,k}(f)\leq c2^{ks}\sum\limits_{j=k+1}^{\infty }2^{\left( j-k\right)
		a(1-\lambda )}(\phi _{j}^{\ast ,a}f)^{1-\lambda }\mathcal{M}(|\mathcal{F}%
	^{-1}\phi _{j}\ast f|^{\lambda }).
	\end{equation*}%
	Using Lemma \ref{lem:lq-inequality}, we obtain that \eqref{second-term} in $%
	\ell ^{\beta }$-quasi-norm can be estimated from above by 
	\begin{align*}
	& c\Big(\sum\limits_{j=0}^{\infty }2^{js\beta }(\phi _{j}^{\ast
		,a}f)^{(1-\lambda )\beta }(\mathcal{M}(|\mathcal{F}^{-1}\phi _{j}\ast
	f|^{\lambda }))^{\beta }\Big)^{1/\beta } \\
	& \lesssim \Big(\sum\limits_{j=0}^{\infty }2^{js\beta }(\phi _{j}^{\ast
		,a}f)^{\beta }\Big)^{(1-\lambda )/\beta }\Big(\sum\limits_{j=0}^{\infty
	}2^{js\beta }(\mathcal{M}(|\mathcal{F}^{-1}\phi _{j}\ast f|^{\lambda
	}))^{\beta /\lambda }\Big)^{\lambda /\beta }.
	\end{align*}%
	Applying the $\dot{K}_{q}^{\alpha ,p}$-quasi-norm and using H\"{o}lder's
	inequality we obtain that 
	\begin{equation*}
	\big\|\Big(\sum\limits_{j=0}^{\infty }(J_{2,k}(f))^{\beta }\Big)^{1/\beta }%
	\big\|_{\dot{K}_{q}^{\alpha ,p}}
	\end{equation*}%
	is bounded by 
	\begin{align*}
	& c\Big\|\Big(\sum\limits_{j=0}^{\infty }2^{js\beta }(\phi _{j}^{\ast
		,a}f)^{\beta }\Big)^{(1-\lambda )/\beta }\Big\|_{\dot{K}_{\frac{q}{1-\lambda 
		}}^{\alpha (1-\lambda ),\frac{p}{1-\lambda }}}\times  \\
	& \Big\|\Big(\sum\limits_{j=0}^{\infty }2^{js\beta }\big(\mathcal{M}(|%
	\mathcal{F}^{-1}\phi _{j}\ast f|^{\lambda })\big)^{\beta /\lambda }\Big)%
	^{\lambda /\beta }\Big\|_{\dot{K}_{\frac{q}{\lambda }}^{\alpha \lambda ,%
			\frac{p}{\lambda }}} \\
	& \lesssim \Big\|\Big(\sum\limits_{j=0}^{\infty }2^{js\beta }(\phi
	_{j}^{\ast ,a}f)^{\beta }\Big)^{1/\beta }\Big\|_{\dot{K}_{q}^{\alpha
			,p}}^{1-\lambda }\Big\|\Big(\sum\limits_{j=0}^{\infty }2^{js\beta }|\mathcal{%
		F}^{-1}\phi _{j}\ast f|^{\beta }\Big)^{1/\beta }\Big\|_{\dot{K}_{q}^{\alpha
			,p}}^{\lambda } \\
	& \lesssim \big\|f\big\|_{\dot{K}_{q}^{\alpha ,p}F_{\beta }^{s}},
	\end{align*}%
	where we have used Lemma \ref{Maximal-Inq} and Theorem \ref{Peetre maximal function}.
	
	$\bullet $\textit{\ Case 2. }$\min (q,\beta )>1$. Assume that $\alpha \geq
	n(1-\frac{1}{q})$. Then we choose $\lambda $ as in \eqref{lamda6}\ and $%
	\frac{n}{q}+\alpha <a<\frac{s}{1-\lambda }$. If $-\frac{n}{q}<\alpha <n(1-%
	\frac{1}{q})$, then we choose $\lambda =1$. The desired estimate can be done
	in the same manner as in Case 1.
	
	\textit{Substep 2.2.} We will estimate%
	\begin{equation*}
	\Big\|\Big(\sum_{k=-\infty }^{-1}2^{sk\beta }|d_{2^{-k}}^{M}f|^{\beta }\Big)%
	^{1/\beta }\Big\|_{\dot{K}_{q}^{\alpha,p }}.
	\end{equation*}%
	We employ the same notations as in Subtep 1.1. Define%
	\begin{equation*}
	H_{k,2}(f)(x)=\int_{B}\big|\sum_{j=0}^{\infty }\Delta _{z2^{-k}}^{M}(\mathcal{F}^{-1}\varphi
	_{j}\ast f)(x)\big|dz,\quad k\leq 0,  x \in \mathbb{R}^{n} .
	\end{equation*}%
	As in the estimation of $J_{2,k}$, we obtain that%
	\begin{equation*}
	H_{2,k}(f)\lesssim 2^{-k\frac{a}{\sigma }(1-\lambda )}\sup_{j\in \mathbb{N}%
		_{0}}\Big(\big(2^{js}(\mathcal{F}^{-1}\varphi _{j})^{\ast ,\frac{a}{\sigma }}f\big)^{1-\lambda
	}\mathcal{M}\big(2^{js}|\mathcal{F}^{-1}\varphi _{j}\ast f|\big)^{\lambda }\Big)
	\end{equation*}%
	and this yields that%
	\begin{equation*}
	\Big(\sum_{k=-\infty }^{-1}2^{sk\beta }|H_{2,k}|^{\beta }\Big)^{1/\beta
	}\lesssim \sup_{j\in \mathbb{N}_{0}}\Big(\big(2^{js}(\mathcal{F}^{-1}\varphi _{j})^{\ast ,%
		\frac{a}{\sigma }}f\big)^{1-\lambda }\mathcal{M}\big(2^{js}|\mathcal{F}^{-1}\varphi _{j}\ast
	f|\big)^{\lambda }\Big).
	\end{equation*}%
	By the same arguments as used in Subtep 2.1 we obtain the desired estimate.  \\
	\textit{Step 3.} Let $f\in \dot{K}_{q}^{\alpha ,p}A_{\beta }^{s}$. We will
	to prove that%
	\begin{equation*}
	\big\|f\big\|_{\dot{K}_{q}^{\alpha ,p}A_{\beta }^{s}}\lesssim \big\|f\big\|_{%
		\dot{K}_{q}^{\alpha ,p}A_{\beta }^{s}}^{\ast }.
	\end{equation*}%
	As the proof for $\dot{K}_{q}^{\alpha ,p}B_{\beta }^{s}$ is similar, we only
	consider $\dot{K}_{q}^{\alpha ,p}F_{\beta }^{s}$. Let $\Psi $ be a function
	in $\mathcal{S}(\mathbb{R}^{n})$ satisfying $\Psi (x)=1$ for $\lvert x\rvert
	\leq 1$ and $\Psi (x)=0$ for $\lvert x\rvert \geq \frac{3}{2}$, and in
	addition radialsymmetric. We make use of an observation made by Nikol'skij 
	\cite{Nikolskii1975} (see also \cite{Si99} and \cite[Section 3.3.2]%
	{Triebel83}). We put%
	\begin{equation*}
	\psi (x)=(-1)^{M+1}\sum\limits_{i=0}^{M-1}(-1)^{i}C_{i}^{M}\Psi (x\left(
	M-i\right) ).
	\end{equation*}%
	The function $\psi $ satisfies $\psi \left( x\right) =1$ for $\left\vert
	x\right\vert \leq \frac{1}{M}$ and $\psi \left( x\right) =0$ for $\left\vert
	x\right\vert \geq \frac{3}{2}$. Then, taking $\varphi _{0}(x)=\psi
	(x),\varphi _{1}(x)=\psi (\frac{x}{2})-\psi (x)$ and $\varphi
	_{j
	}(x)=\varphi _{1}(2^{-j+1}x)$ for $j=2,3,...$, we obtain that $\{
	\varphi _{j}\}_{j\in \mathbb{N}_{0}} $\ is a smooth dyadic resolution of unity. This yields
	that%
	\begin{equation*}
	\Big\|\Big(\sum_{j=0}^{\infty }2^{js\beta }|\mathcal{F}^{-1}\varphi _{j}\ast
	f|^{\beta }\Big)^{1/\beta }\Big\|_{\dot{K}_{q}^{\alpha ,p}}
	\end{equation*}%
	is a quasi-norm equivalent in $\dot{K}_{q}^{\alpha ,p}{F_{\beta }^{s}}$. Let
	us prove that the last expression is bounded by%
	\begin{equation}
	C\big\Vert f\big\Vert_{\dot{K}_{q}^{\alpha ,p}{F_{\beta }^{s}}}^{\ast }.
	\label{second-est}
	\end{equation}%
	We observe that%
	\begin{equation*}
	\mathcal{F}^{-1}\varphi _{0}\ast f(x)=\left( -1\right) ^{M+1}\int_{\mathbb{R}%
		^{n}}\mathcal{F}^{-1}\Psi \left( z\right) \Delta _{-z}^{M}f(x)dz+f(x)\int_{%
		\mathbb{R}^{n}}\mathcal{F}^{-1}\Psi \left( z\right) dz
	\end{equation*}%
	Moreover, it holds for $x\in \mathbb{R}^{n}$ and $j=1,2,...$%
	\begin{equation*}
	\mathcal{F}^{-1}\varphi _{j}\ast f\left( x\right) =\left( -1\right)
	^{M+1}\int_{\mathbb{R}^{n}}\Delta _{2^{-j}y}^{M}f\left( x\right) \widetilde{%
		\Psi }\left( y\right) dy,
	\end{equation*}%
	with $\widetilde{\Psi } =\mathcal{F}^{-1}\Psi-2^{-n}\mathcal{F}^{-1}\Psi ( \cdot /2) $. Now, for 
	$j\in \mathbb{N}_{0}$ we have%
	\begin{align}
	&\int_{\mathbb{R}^{n}}|\Delta _{2^{-j}y}^{M}f( x) ||\widetilde{%
		\Psi }( y) |dy  \notag \\
	&=\int_{\left\vert y\right\vert \leq 1}|\Delta _{2^{-j}y}^{M}f(
	x) ||\widetilde{\Psi }( y) |dy+\int_{\left\vert
		y\right\vert >1}|\Delta _{2^{-j}y}^{M}f( x) ||\widetilde{\Psi }%
	( y) |dy.  \label{diff2.1}
	\end{align}%
	Thus, we need only to estimate the second term of \eqref{diff2.1}. We write 
	\begin{align}
	&2^{sj}\int_{\left\vert y\right\vert >1}|\Delta _{2^{-j}y}^{M}f(
	x) ||\widetilde{\Psi }( y) |dy  \notag \\
	&=2^{sj}\sum\limits_{k=0}^{\infty }\int_{2^{k}<\left\vert y\right\vert \leq
		2^{k+1}}|\Delta _{2^{-j}y}^{M}f( x) ||\widetilde{\Psi }(
	y) |dv  \notag \\
	&\leq c2^{sj}\sum\limits_{k=0}^{\infty }2^{nj-Nk}\int_{2^{k-j}<\left\vert
		h\right\vert \leq 2^{k-j+1}}|\Delta _{h}^{M}f( x) |dh
	\label{diff}
	\end{align}%
	where $N>0$ is at our disposal and we have used the properties of the
	function $\widetilde{\Psi }$, $|\widetilde{\Psi }( x) |\leq
	c(1+\left\vert x\right\vert )^{-N},$ for any $x\in \mathbb{R}^{n}$ and any $%
	N>0$. Without lost of generality, we may assume $1\leq \beta \leq \infty $.
	Now, the right-hand side of \eqref{diff} in $\ell^{\beta }$-norm is bounded
	by 
	\begin{equation}
	c\sum\limits_{k=0}^{\infty }2^{-Nk}\Big(\sum\limits_{j=0}^{\infty
	}2^{(s+n)j\beta }\Big(\int_{\left\vert h\right\vert \leq 2^{k-j+1}}|\Delta
	_{h}^{M}f(x)|dh\Big)^{\beta }\Big)^{1/\beta }.  \label{diff1}
	\end{equation}%
	After a change of variable $j-k-1=v$, we estimate \eqref{diff1} \ by%
	\begin{equation*}
	c\sum\limits_{k=0}^{\infty }2^{(s+n-N)k}\Big(\sum\limits_{v=-k-1}^{\infty
	}2^{sv\beta }\big(d_{2^{-v}}^{M}f(x)\big)^{\beta }\Big)^{1/\beta }\lesssim %
	\Big(\sum\limits_{v=-{\infty }}^{\infty
	}2^{sv\beta }\big(d_{2^{-v}}^{M}f(x)\big)^{\beta }\Big)^{1/\beta },
	\end{equation*}%
	where we choose $N>n+s$. Taking the  $\dot{K}_{q}^{\alpha ,p}$-quasi-norm we obtain the desired estimate \eqref{second-est}. The proof is complete.
\end{proof}

We would like to mention that%
\begin{equation}
\left\Vert f(\lambda \cdot )\right\Vert _{\dot{K}_{q}^{\alpha ,p}B_{\beta
	}^{s}}^{\ast }\approx\lambda ^{-\alpha -\frac{n}{q}}\big\|f\big\|_{\dot{K}%
	_{q}^{\alpha ,p}}+\lambda ^{s-\alpha -\frac{n}{q}}\Big(\int_{0}^{\infty
}t^{-s\beta }\big\|d_{t}^{M}f\big\|_{\dot{K}_{q}^{\alpha ,p}}^{\beta }\frac{%
	dt}{t}\Big)^{\frac{1}{\beta }}  \label{delation1}
\end{equation}%
and%
\begin{equation*}
\left\Vert f(\lambda \cdot )\right\Vert _{\dot{K}_{q}^{\alpha ,p}F_{\beta
	}^{s}}^{\ast }\approx\lambda ^{-\alpha -\frac{n}{q}}\big\|f\big\|_{\dot{K}%
	_{q}^{\alpha ,p}}+\lambda ^{s-\alpha -\frac{n}{q}}\big\|\Big(%
\int_{0}^{\infty }t^{-s\beta }(d_{t}^{M}f)^{\beta }\frac{dt}{t}\Big)^{\frac{1%
	}{\beta }}\big\|_{\dot{K}_{q}^{\alpha ,p}}
\end{equation*}%
for any $\lambda >0,0<p\leq \infty ,0<q\leq \infty ,\alpha >-\frac{n}{q}%
,\max (\sigma _{q},\alpha -\alpha _{0})<s<M\ $($0<p,q<\infty $ and $\max
(\sigma _{q,\beta },\alpha -\alpha _{0})<s<M$ in the $\dot{K}F$-case) and $%
M\in \mathbb{N}$. \\

Let $\varphi ^{j}(x)=\varphi _{0}(2^{-j}x)-\varphi _{0}(2^{1-j}x)$ for $j\in 
\mathbb{Z}$ and $x\in \mathbb{R}^{n}$. In view of \cite{XuYang03} we have
the following equivalent norm of $\dot{K}_{q}^{\alpha ,p}$. Let\textit{\ }$%
1<p,q<\infty \ $and $-\frac{n}{q}<\alpha <n-\frac{n}{q}$\textit{. }Then%
\begin{equation}
\Big\|\Big(\sum\limits_{j=-\infty }^{\infty }\left\vert \mathcal{F}%
^{-1}\varphi ^{j}\ast f\right\vert ^{2}\Big)^{1/2}\Big\|_{\dot{K}%
	_{q}^{\alpha ,p}}\approx \big\|f\big\|_{\dot{K}_{q}^{\alpha ,p}},
\label{herz-norm}
\end{equation}%
holds for all $f\in\dot{K}_{q}^{\alpha ,p}$.

Let $s\in\mathbb{R},0<p,q<\infty ,0<\beta \leq \infty \ $ and $\alpha >-\frac{n}{q}$. We set
\begin{equation*}
\big\|f\big\|_{\dot{K}_{q}^{\alpha ,p}\dot{B}_{\beta }^{s}}=\Big(%
\sum\limits_{j=-\infty }^{\infty }2^{js\beta }\big\|\mathcal{F}^{-1}\varphi
^{j}\ast f\big\|_{\dot{K}_{q}^{\alpha ,p}}^{\beta }\Big)^{1/\beta }
\end{equation*}
and
\begin{equation*}
\big\|f\big\|_{\dot{K}_{q}^{\alpha ,p}\dot{F}_{\beta }^{s}}=\Big\|\Big(%
\sum\limits_{j=-\infty }^{\infty }2^{js\beta }\left\vert \mathcal{F}%
^{-1}\varphi ^{j}\ast f\right\vert ^{\beta }\Big)^{1/\beta }\Big\|_{\dot{K}%
	_{q}^{\alpha ,p}}.
\end{equation*}

\begin{proposition}
	\label{Triebel1 copy(3)} Let\textit{\ }$s>\max (\sigma _{q},\alpha -n+\frac{n%
	}{q}),0<p,q<\infty ,0<\beta \leq \infty \ $and $\alpha >-\frac{n}{q}$\textit{%
		. }\newline
	$\mathrm{(i)}$ Let $s>\max (\sigma _{q},\alpha -n+\frac{n%
	}{q})$ and  $f\in \dot{K}_{q}^{\alpha ,p}B_{\beta }^{s}$. \textit{Then%
	}%
	\begin{equation*}
	\big\|f\big\|_{\dot{K}_{q}^{\alpha ,p}B_{\beta }^{s}}\approx \big\|f\big\|_{%
		\dot{K}_{q}^{\alpha ,p}}+\big\|f\big\|_{\dot{K}_{q}^{\alpha ,p}\dot{B}%
		_{\beta }^{s}},
	\end{equation*}%
	$\mathrm{(ii)}$ \textit{Let } $s>\max (\sigma _{q,\beta},\alpha -n+\frac{n%
	}{q})$ and $f\in \dot{K}_{q}^{\alpha ,p}F_{\beta }^{s}$. 
	\textit{Then}%
	\begin{equation*}
	\big\|f\big\|_{\dot{K}_{q}^{\alpha ,p}F_{\beta }^{s}}\approx \big\|f\big\|_{%
		\dot{K}_{q}^{\alpha ,p}}+\big\|f\big\|_{\dot{K}_{q}^{\alpha ,p}\dot{F}%
		_{\beta }^{s}}.
	\end{equation*}%
	
\end{proposition}
\begin{proof}
	As the proof for (i) is similar, we only consider  (ii). We use the following Marschall's inequality which given in \cite[Proposition 1.5]{Mar91}, see also \cite{DH19}. Let\textit{\ }$A>0,R\geq 1$. Let $b\in \mathcal{D}(%
	\mathbb{R}
	^{n})$ and a function $g\in C^{\infty }(%
	\mathbb{R}
	^{n})$ such that 
	\begin{equation*}
	\mathrm{supp}\mathcal{F}g\subseteq \left\{ \xi \in \mathbb{R}^{n}:\left\vert
	\xi \right\vert \leq AR\right\} \quad \text{and}\quad \mathrm{supp}%
	b\subseteq \left\{ \xi \in \mathbb{R}^{n}:\left\vert \xi \right\vert \leq
	A\right\} .
	\end{equation*}%
	Then%
	\begin{equation*}
	\left\vert \mathcal{F}^{-1}b\ast g(x)\right\vert \leq c (AR)^{\frac{n}{t}-n}\left\Vert b\right\Vert _{\dot {B}_{1,t}^{\frac{n}{t}}}\mathcal{M}_{t}(g)(x)
	\end{equation*}%
	for any $0<t\leq 1$ and any $x\in \mathbb{R}^{n}$, where $c$ is independent
	of $A$, $R$, $x$, $b$, $j$ and $g$. Here $\dot {B}_{1,t}^{\frac{n}{t}} $ denotes the homogeneous Besov spaces. We have%
	\begin{equation*}
	\mathcal{F}^{-1}\varphi ^{j}\ast f=\mathcal{F}^{-1}\varphi ^{j}\ast \mathcal{%
		F}^{-1}\varphi _{0}\ast f,\quad -j\in \mathbb{N}.
	\end{equation*}%
	Therefore,%
	\begin{equation*}
	\left\vert \mathcal{F}^{-1}\varphi ^{j}\ast f(x)\right\vert \leq c\left\Vert \varphi ^{j}\right\Vert _{\dot {B}_{1,t}^{\frac{n}{t}}}\mathcal{M}_{t}(\mathcal{F}^{-1}\varphi
	_{0}\ast f)(x) 
	\leq c 2^{j(n-\frac{n}{t})}\mathcal{M}_{t}(\mathcal{F}^{-1}\varphi _{0}\ast f)(x),\quad x\in 
	\mathbb{R}^{n},
	\end{equation*}%
	where the positive constant $c$ is independent of $j$ and $x$. If we choose $\frac{n}{s +n}<t<\min (1,q,\beta,\frac{n}{\alpha +\frac{n}{q}})$ then%
	\begin{equation*}
	\Big(\sum\limits_{j=-\infty }^{-1}2^{js\beta }\left\vert \mathcal{F}%
	^{-1}\varphi ^{j}\ast f\right\vert ^{\beta }\Big)^{1/\beta }\lesssim 
	\mathcal{M}_{t}(\mathcal{F}^{-1}\varphi _{0}\ast f).
	\end{equation*}%
	Taking the $\dot{K}_{q}^{\alpha ,p}$-quasi-norm$\ $and using %
	\eqref{convolution estimate} we obtain 
	\begin{equation*}
	\Big\|\Big(\sum\limits_{j=-\infty }^{\infty }2^{js\beta }\left\vert \mathcal{%
		F}^{-1}\varphi ^{j}\ast f\right\vert ^{\beta }\Big)^{1/\beta }\Big\|_{\dot{K}%
		_{q}^{\alpha ,p}}\lesssim \big\|f\big\|_{\dot{K}_{q}^{\alpha ,p}F_{\beta
		}^{s}}.
	\end{equation*}%
	Because of $s>\max (\sigma _{q},\alpha -n+\frac{n}{q})$ the series $%
	\sum\limits_{j=0}^{\infty }\mathcal{F}^{-1}\varphi _{j}\ast f$ converges not
	only in $\mathcal{S}^{\prime }(\mathbb{R}^{n})$ but almost everywhere in $%
	\mathbb{R}^{n}$. Then%
	\begin{equation*}
	\big\|f\big\|_{\dot{K}_{q}^{\alpha ,p}}\lesssim \big\|\mathcal{F}%
	^{-1}\varphi _{0}\ast f\big\|_{\dot{K}_{q}^{\alpha ,p}}+\Big(%
	\sum\limits_{j=1}^{\infty }\big\|\mathcal{F}^{-1}\varphi _{j}\ast f\big\|_{%
		\dot{K}_{q}^{\alpha ,p}}^{\min (1,p,q)}\Big)^{1/\min (1,p,q)}.
	\end{equation*}%
	Therefore $\big\|f\big\|_{\dot{K}_{q}^{\alpha ,p}}+\big\|f\big\|_{\dot{K}%
		_{q}^{\alpha ,p}\dot{F}_{\beta }^{s}}$ can be estimated from above by $c%
	\big\|f\big\|_{\dot{K}_{q}^{\alpha ,p}F_{\beta }^{s}}$. Obviously%
	\begin{equation*}
	\mathcal{F}^{-1}\varphi _{0}\ast f =\sum\limits_{j=0}^{N}\mathcal{F}%
	^{-1}\varphi _{j}\ast f-\sum\limits_{j=1}^{N}\mathcal{F}^{-1}\varphi
	_{j}\ast f 
	=g_{N}+h_{N},\quad N\in \mathbb{N}.
	\end{equation*}%
	We have%
	\begin{equation*}
	\big\|h_{N}\big\|_{\dot{K}_{q}^{\alpha ,p}}\leq \Big(\sum\limits_{j=1}^{%
		\infty }\big\|\mathcal{F}^{-1}\varphi ^{j}\ast f\big\|_{\dot{K}_{q}^{\alpha
			,p}}^{\min (1,p,q)}\Big)^{1/\min (1,p,q)},\quad N\in \mathbb{N}.
	\end{equation*}%
	By Lebesgue's dominated convergence theorem, it follows that\ $\big\|g_{N}-f%
	\big\|_{\dot{K}_{q}^{\alpha ,p}}$\ tends to zero as $N$ tends to infinity.
	Therefore $\big\|\mathcal{F}^{-1}\varphi _{0}\ast f\big\|_{\dot{K}%
		_{q}^{\alpha ,p}}$ can be estimated from above by the quasi-norm 
	\begin{equation*}
	c\big\|f\big\|_{\dot{K}_{q}^{\alpha ,p}}+c\big\|f\big\|_{\dot{K}_{q}^{\alpha
			,p}\dot{F}_{\beta }^{s}}.
	\end{equation*}%
	The proof is complete.
\end{proof}

\begin{proposition}
	\label{Triebel1 copy(1)} Let\textit{\ }$s>0,1<p,q<\infty \ $and $-\frac{n}{q}%
	<\alpha <n-\frac{n}{q}$\textit{. Let}%
	\begin{equation*}
	\mathcal{S}_{0}(\mathbb{R}^{n})=\big\{f\in \mathcal{S}(\mathbb{R}^{n}):%
	\mathrm{supp}\mathcal{F}f\in \mathbb{R}^{n}\backslash \{0\}\big\}.
	\end{equation*}%
	Then $\mathcal{S}_{0}(\mathbb{R}^{n})$ is dense in $\dot{k}_{q,s}^{\alpha
		,p} $.
\end{proposition}
\begin{proof}
	Let $\varphi _{0}=\varphi $\ be as above. As in \cite%
	{Triebel13} it is sufficient to approximate $f\in \mathcal{S}(\mathbb{R}%
	^{n}) $ in $\dot{W}_{q,k}^{\alpha ,p}$, $k\in \mathbb{N}$, by functions
	belonging to $\mathcal{S}_{0}(\mathbb{R}^{n})$. We have%
	\begin{equation*}
	|D^{\alpha }\mathcal{F}^{-1}(\varphi (2^{j}\cdot )\mathcal{F}f)|=2^{-jn}|%
	\tilde{\varphi}_{j}\ast D^{\alpha }f|\leq 2^{-jn}\mathcal{M}(\tilde{\varphi}%
	_{j}),
	\end{equation*}%
	where $\tilde{\varphi}_{j}=\mathcal{F}^{-1}\varphi (2^{-j}\cdot ),j\in 
	\mathbb{N}\ $and $\alpha \in \mathbb{N}^{n}$. From 
	\eqref{convolution estimate} we obtain 
	\begin{equation*}
	\big\|D^{\alpha }\mathcal{F}^{-1}(\varphi (2^{j}\cdot )\mathcal{F}f)\big\|_{%
		\dot{K}_{q}^{\alpha ,p}} \leq c2^{-jn}\big\|\tilde{\varphi}_{j}\big\|_{%
		\dot{K}_{q}^{\alpha ,p}} \\
	\leq c2^{j(\frac{n}{q}-n+\alpha )},
	\end{equation*}%
	where the positive constant $c$ is independent of $j$. Since $\alpha <n-%
	\frac{n}{q}$, we obtain that $f-\mathcal{F}^{-1}(\varphi (2^{j}\cdot )%
	\mathcal{F}f)$\ approximate $f\in \mathcal{S}(\mathbb{R}^{n})$ in $\dot{W}%
	_{q,k}^{\alpha ,p}$, $k\in \mathbb{N}$. The proof of the proposition is
	complete.
\end{proof}

\begin{proposition}
	\label{Triebel1 copy(2)} Let\textit{\ }$s>0,1<p,q<\infty \ $and $-\frac{n}{q}%
	<\alpha <n-\frac{n}{q}$\textit{. Let }$f\in \dot{k}_{q,s}^{\alpha ,p}$. 
	\textit{Then}%
	\begin{equation*}
	\big\|f\big\|_{\dot{k}_{q,s}^{\alpha ,p}}\approx \big\|f\big\|_{\dot{K}%
		_{q}^{\alpha ,p}}+\big\|(-\Delta )^{\frac{s}{2}}f\big\|_{\dot{K}_{q}^{\alpha
			,p}},
	\end{equation*}%
	where%
	\begin{equation*}
	(-\Delta )^{\frac{s}{2}}f=\mathcal{F}^{-1}(|\xi |^{s}\mathcal{F}f).
	\end{equation*}
\end{proposition}
\begin{proof}
	Let $f\in \mathcal{S}_{0}(\mathbb{R}^{n})$. We apply
	Marschall's inequality to $g_{j}=\mathcal{F}^{-1}(\varphi ^{j}|x|^{s}%
	\mathcal{F}f),j\in \mathbb{Z}$ and $b_{j}(x)=2^{js}|x|^{-s}\psi ^{j}(x),j\in 
	\mathbb{Z},x\in \mathbb{R}^{n}$ where%
	\begin{equation*}
	\varphi ^{j}(x)=\varphi _{0}(2^{-j}x)-\varphi _{0}(2^{1-j}x),\quad \psi
	^{j}=\varphi ^{j-1}+\varphi ^{j}+\varphi ^{j+1},\quad j\in \mathbb{Z},x\in 
	\mathbb{R}^{n}.
	\end{equation*}%
	Then%
	\begin{equation*}
	\left\vert \mathcal{F}^{-1}b_{j}\ast g_{j}(x)\right\vert \leq c\left\Vert
	b_{j}\right\Vert _{B_{1,1}^{n}}\mathcal{M}(%
	\mathcal{F}^{-1}(\varphi ^{j}|\xi |^{s}\mathcal{F}f))(x) 
	\leq c\mathcal{M}(\mathcal{F}^{-1}(\varphi ^{j}|\xi |^{s}\mathcal{F}f))(x)
	\end{equation*}%
	for any $j\in \mathbb{Z}$ and any $x\in \mathbb{R}^{n}$, where $c$ is
	independent of $j$. Let $j\in \mathbb{Z} $. In view of the fact that%
	\begin{equation*}
	\mathcal{F}^{-1}\varphi ^{j}\ast f=\mathcal{F}^{-1}(\varphi ^{j}\mathcal{F}%
	f)=2^{-js}\mathcal{F}^{-1}(2^{js}|\xi |^{-s}\psi ^{j}|x|^{s}\varphi ^{j}%
	\mathcal{F}f)=2^{-js}\mathcal{F}^{-1}(b_{j}|\xi |^{s}\varphi ^{j}\mathcal{F}%
	f),
	\end{equation*}%
	by Lemma \ref{Maximal-Inq} and \eqref{herz-norm} we obtain%
	\begin{align*}
	\Big\|\Big(\sum\limits_{j=-\infty }^{\infty }2^{2sj}\left\vert \mathcal{F}%
	^{-1}\varphi ^{j}\ast f\right\vert ^{2}\Big)^{1/2}\Big\|_{\dot{K}%
		_{q}^{\alpha ,p}} &\lesssim \Big\|\Big(\sum\limits_{j=-\infty }^{\infty
	}\left\vert \mathcal{F}^{-1}(\varphi ^{j}|\xi |^{s}\mathcal{F}f)\right\vert
	^{2}\Big)^{1/2}\Big\|_{\dot{K}_{q}^{\alpha ,p}} \\
	&\lesssim \big\|\mathcal{F}^{-1}(|\xi |^{s}\mathcal{F}f)\big\|_{\dot{K}%
		_{q}^{\alpha ,p}}.
	\end{align*}%
	The same arguments can be used to prove the opposite inequality in view of
	the fact that%
	\begin{equation*}
	\mathcal{F}^{-1}(\varphi ^{j}|\xi |^{s}\mathcal{F}f)=\mathcal{F}%
	^{-1}(2^{-js}\psi ^{j}|\xi |^{s}2^{js}\varphi ^{j}\mathcal{F}f)=\mathcal{F}%
	^{-1}(b_{j}2^{js}\varphi ^{j}\mathcal{F}f),\quad j\in \mathbb{Z}.
	\end{equation*}%
	The rest follows by Propositions \ref{Triebel1 copy(3)} and \ref{Triebel1
		copy(1)}. The proof is complete. 
\end{proof}
\begin{definition}
	{\rm	Let $0<u\leq p<\infty $. The Morrey space $M_{u}^{p}$ is defined
		to be the set of all $u$-locally Lebesgue-integrable functions $f$ on $%
		\mathbb{R}^{n}$ such that%
		\begin{equation*}
		\left\Vert f\right\Vert _{M_{u}^{p}}=\sup \left\vert B\right\vert
		^{\frac{1}{p}-\frac{1}{u}}\big\|f\chi _{B}\big\|_{u}<\infty ,
		\end{equation*}%
		where the supremum is taken over all balls $B$ in $\mathbb{R}^{n}$.}
\end{definition}

\begin{remark}
	{\rm	The Morrey spaces $M_{u}^{p}$ are quasi-Banach spaces, Banach
		spaces for $u\geq 1$, were introduced by Morrey to study some PDE's, see 
		\cite{Mo38}. One can easily seen that  $M_{p}^{p}=L^{p}$ and that for $0<u\leq v\leq p<\infty $,%
		\begin{equation*}
		M_{v}^{p}\hookrightarrow M_{u}^{p}.
		\end{equation*}}
\end{remark}
The Sobolev Morrey spaces are defined as follows.

\begin{definition}
	{\rm	Let $1<u\leq p<\infty $ and $m=1,2,...$. The Sobolev Morrey space $M_{u}^{m,p}$ is defined to be the set of all $u$-locally Lebesgue-integrable
		functions $f$ on $\mathbb{R}^{n}$ such that%
		\begin{equation*}
		\left\Vert f\right\Vert _{M_{u}^{m,p}}=\left\Vert f\right\Vert _{%
			M_{u}^{p}}+\sum\limits_{|\alpha |\leq m}\left\Vert D^{\alpha
		}f\right\Vert _{M_{u}^{p}}<\infty .
		\end{equation*}}
\end{definition}

Let now recall the definition of Besov-Morrey and Triebel-Lizorkin-Morrey\
spaces. Let $\{\varphi _{j}\}_{j\in \mathbb{N}_{0}}$ be a resolution of unity, see Section 2.

\begin{definition}
	{\rm	Let $s\in \mathbb{R},0<u\leq p<\infty $\ and $0<q\leq \infty $.The
		Besov-Morrey\ space\ $\mathcal{N}_{p,q,u}^{s}$\ is the collection of all $%
		f\in \mathcal{S}^{\prime }(\mathbb{R}^{n})$\ such that 
		\begin{equation*}
		\left\Vert f\right\Vert _{\mathcal{N}_{p,q,u}^{s}}=\Big(\sum\limits_{j=0}^{%
			\infty }2^{jsq}\big\|\mathcal{F}^{-1}\varphi _{j}\ast f\big\|_{M_{u}^{p}}^{q}\Big)^{1/q}<\infty .
		\end{equation*}%
		In the limiting case $q=\infty $ the usual modification is required.\textit{%
			\newline
		}The Triebel-Lizorkin-Morrey space $\mathcal{E}_{p,q,u}^{s}$\ is the
		collection of all $f\in \mathcal{S}^{\prime }(\mathbb{R}^{n})$\ such that 
		\begin{equation*}
		\left\Vert f\right\Vert _{\mathcal{E}_{p,q,u}^{s}}=\Big\|\Big(%
		\sum\limits_{j=0}^{\infty }2^{jsq}\left\vert \mathcal{F}^{-1}\varphi
		_{j}\ast f\right\vert ^{q}\Big)^{1/q}\Big\|_{M_{u}^{p}}<\infty .
		\end{equation*}%
		In the limiting case $q=\infty $ the usual modification is required.}
\end{definition}
We have%
\begin{equation*}
\mathcal{E}_{p,2,u}^{m}=M_{u}^{m,p},\quad m\in \mathbb{N},\quad
1<u\leq p<\infty ,  \label{coincidence}
\end{equation*}%
with equivalent norms, see \cite[Theorem 3.1]{Si12}. In particular, we have that%
\begin{equation}
\mathcal{E}_{p,2,u}^{0}=M_{u}^{p},\quad 1<u\leq p<\infty, \label{Morrey}
\end{equation}%
also in the sense of with equivalent norms, see \cite[Proposition 4.1]{Ma06}.

\begin{theorem}
	\label{Sobolev-embeddings}Let $s_{i}\in \mathbb{R},0<q_{i}\leq \infty
	,0<u_{i}\leq p_{i}<\infty ,i=1,2$. There is a continuous embedding%
	\begin{equation*}
	\mathcal{E}_{p_{1},q_{1},u_{1}}^{s_{1}}\hookrightarrow \mathcal{E}%
	_{p_{2},q_{2},u_{2}}^{s_{2}}
	\end{equation*}%
	if, and only if,%
	\begin{equation*}
	p_{1}\leq p_{2}\quad \text{and}\quad \frac{u_{2}}{p_{2}}\leq \frac{u_{1}}{%
		p_{1}}
	\end{equation*}%
	and%
	\begin{equation*}
	s_{1}-\frac{n}{p_{1}}>s_{2}-\frac{n}{p_{2}}\quad or\quad s_{1}-\frac{n}{p_{1}}=s_{2}-\frac{n}{p_{2}}\quad \text{and}\quad p_{1}\neq
	p_{2}.
	\end{equation*}%
	
\end{theorem}

For the proof of these Sobolev embeddings, see \cite[Theorem 3.1]{HS14}.

\begin{remark}
	{\rm A detailed study of Besov-Morrey and Triebel-Lizorkin-Morrey
		spaces including their history and properties can be found in \cite{HS14, Ma06, Ma05, Si12, SiYY} and references therein.}
	
\end{remark}

\section{ \large  Caffarelli-Kohn-Nirenberg  inequalities}

As mentioned in the introduction, Caffarelli-Kohn-Nirenberg\ inequalities
play a crucial role to study regularity and integrability for solutions of
nonlinear partial differential equations, see \cite{FS08, Xuan05}.
The main aim of this section is to extend these inequalities to more general
function spaces. Let $\{\varphi _{j}\}_{j\in \mathbb{N}_{0}}$\ be a
resolution of unity and%
\begin{equation*}
Q_{J}f=\sum\limits_{j=0}^{J}\mathcal{F}^{-1}\varphi _{j}\ast f,\quad J\in 
\mathbb{N},f\in \mathcal{S}^{\prime }(\mathbb{R}^{n}).
\end{equation*}

\subsection{\large  CKN inequalities in Herz-type Besov and Triebel-Lizorkin spaces}

In this section, we investigate the Caffarelli, Kohn and Nirenberg
inequalities in $\dot{K}_{q}^{\alpha ,p}A_{\beta }^{s}$ spaces. The main
results of this section based on the following proposition.

\begin{proposition}
	\label{Triebel1}\textit{Let }$\alpha _{1},\alpha _{2}\in \mathbb{R},\sigma
	\geq 0,1<r,v<\infty ,0<\tau ,u\leq \infty $ \textit{and} 
	\begin{equation*}
	-\frac{n}{v}<\alpha _{1}<n-\frac{n}{v}.
	\end{equation*}%
	$\mathrm{(i)}$ \textit{Assume that }$1<u\leq v<\infty $ and $\alpha _{2}\geq
	\alpha _{1}$. \textit{Then for all }$f\in \dot{K}_{u}^{\alpha _{2},\delta
	}\cap \mathcal{S}^{\prime }(\mathbb{R}^{n})$ and all $J\in \mathbb{N}$,%
	\begin{equation}
	\big\|Q_{J}f\big\|_{\dot{k}_{v,\sigma }^{\alpha _{1},r}}\leq c2^{J(\frac{n}{u%
		}-\frac{n}{v}+\alpha _{2}-\alpha _{1}+\sigma )}\big\|f\big\|_{\dot{K}%
		_{u}^{\alpha _{2},\delta }},  \label{key-estimate}
	\end{equation}%
	where%
	\begin{equation*}
	\delta =\left\{ 
	\begin{array}{ccc}
	r, & \text{if} & \alpha _{2}=\alpha _{1}, \\ 
	\tau , & \text{if} & \alpha _{2}>\alpha _{1}%
	\end{array}%
	\right.  \label{condition2}
	\end{equation*}%
	\textit{and the positive constant }$c$\textit{\ is independent of }$J$%
	\textit{.}$\newline
	\mathrm{(ii)}$ \textit{Assume that }$1<v\leq u<\infty $ and $\alpha _{2}\geq
	\alpha _{1}+\frac{n}{v}-\frac{n}{u}$. \textit{Then for all }$f\in \dot{K}%
	_{u}^{\alpha _{2},\delta }\cap \mathcal{S}^{\prime }(\mathbb{R}^{n})$ and all $J\in \mathbb{N}$, \eqref{key-estimate} holds where \textit{the positive constant }$c$\textit{\ is independent of }$J$ and%
	\begin{equation*}
	\delta =\left\{ 
	\begin{array}{ccc}
	r, & \text{if} & \alpha _{2}=\alpha _{1}+\frac{n}{v}-\frac{n}{u}, \\ 
	\tau , & \text{if} & \alpha _{2}>\alpha _{1}+\frac{n}{v}-\frac{n}{u}.%
	\end{array}%
	\right.  \label{condition3}
	\end{equation*}%
	
\end{proposition}

\begin{proof} By similarity, we only give the proof for
	(i). Let $\sigma =\theta m+(1-\theta )0$, $\alpha \in \mathbb{N}^{n}$ with $%
	0<\theta <1$ and $|\alpha |\leq m$. From \eqref{Interpolation} we have%
	\begin{equation*}
	\big\|Q_{J}f\big\|_{\dot{K}_{v}^{\alpha _{1},r}A_{2}^{\sigma }}\leq \big\|%
	Q_{J}f\big\|_{\dot{K}_{v}^{\alpha _{1},r}A_{2}^{0}}^{1-\theta }\big\|Q_{J}f%
	\big\|_{\dot{K}_{v}^{\alpha _{1},r}A_{2}^{m}}^{\theta }.
	\end{equation*}%
	Observe that%
	\begin{equation*}
	\dot{K}_{v}^{\alpha _{1},r}A_{2}^{\sigma }=\dot{k}_{v,\sigma }^{\alpha
		_{1},r},\quad \dot{K}_{v}^{\alpha _{1},r}A_{2}^{m}=\dot{W}_{v,m}^{\alpha
		_{1},r},\quad \text{and}\quad \dot{K}_{v}^{\alpha _{1},r}A_{2}^{0}=\dot{K}%
	_{v}^{\alpha _{1},r},
	\end{equation*}%
	see \eqref{coincidence1}, \eqref{coincidence2} and \eqref{coincidence3}. It
	follows that%
	\begin{equation*}
	\big\|Q_{J}f\big\|_{\dot{k}_{v,\sigma }^{\alpha _{1},r}}\leq \big\|Q_{J}f%
	\big\|_{\dot{K}_{v}^{\alpha _{1},r}}^{1-\theta }\big\|Q_{J}f\big\|_{\dot{W}%
		_{v,m}^{\alpha _{1},r}}^{\theta },
	\end{equation*}%
	where the positive constant $c$\ is independent of $J$. Observe that%
	\begin{equation*}
	Q_{J}f=2^{Jn}\mathcal{F}^{-1}\varphi _{0}(2^{J}\cdot )\ast f.
	\end{equation*}%
	Therefore,%
	\begin{equation*}
	D^{\alpha }(Q_{J}f)=2^{J(|\alpha |+n)}\mathcal{\omega }_{J}\ast
	f=2^{J|\alpha |}\tilde{Q}_{J}f,\quad|\alpha |\leq m
	\end{equation*}%
	with $\mathcal{\omega }_{J}(x)=D^{\alpha }(\mathcal{F}^{-1}\varphi
	_{0})(2^{J}x)$, $x\in \mathbb{R}^{n}$. Recall that%
	\begin{equation*}
	|\tilde{Q}_{J}f|\lesssim \mathcal{M}(f).  \label{maximale}
	\end{equation*}%
	Applying Lemma \ref{Bernstein-Herz-ine1} and the estimate 
	\eqref{convolution estimate}, we obtain%
	\begin{align*}
	\big\|D^{\alpha }(Q_{J}f)\big\|_{\dot{K}_{v}^{\alpha _{1},r}} &\leq c2^{J(%
		\frac{n}{u}-\frac{n}{v}+\alpha _{2}-\alpha _{1}+|\alpha |)}\big\|\tilde{Q}%
	_{J}f\big\|_{\dot{K}_{u}^{\alpha _{2},\delta }} \\
	&\leq c2^{J(\frac{n}{u}-\frac{n}{v}+\alpha _{2}-\alpha _{1}+m)}\big\|f\big\|%
	_{\dot{K}_{u}^{\alpha _{2},\delta }}
	\end{align*}%
	for any $|\alpha |\leq m$.
	This finish the proof.
\end{proof}

\begin{remark}
	{\rm	With $\alpha _{1}=\alpha _{2}=0$ and $r=v$\ the  estimation \eqref{key-estimate} can be rewritten as%
		\begin{align*}
		\big\|Q_{J}f\big\|_{H_{v}^{\sigma }} &\leq c2^{J(\frac{n}{u}-\frac{n}{v}%
			+\sigma )}\big\|f\big\|_{\dot{K}_{u}^{0,v}} \\
		&\leq c2^{J(\frac{n}{u}-\frac{n}{v}+\sigma )}\big\|f\big\|_{u},
		\end{align*}%
		because of $1<u\leq v<\infty $ which has been proved by  Triebel in \cite[Proposition 4.5]{Triebel13}.}
\end{remark}

Now we are in position to state the main results of this section.

\begin{theorem}
	\label{Triebel2}Let $0<p,\tau ,\beta ,\varrho < \infty $, $1<r,v,u<\infty ,\sigma \geq 0,$%
	\begin{equation}
	-\frac{n}{v}<\alpha _{1}<n-\frac{n}{v},\quad -\frac{n}{u}<\alpha _{2}<n-%
	\frac{n}{u},\quad \alpha _{3}>-\frac{n}{p},\quad v\geq \max (p,u),
	\label{cond1.1}
	\end{equation}%
	\begin{equation}
	s-\frac{n}{p}+\frac{n}{u}+\alpha _{2}-\alpha _{3}>\sigma -\frac{n}{v}+\alpha
	_{2}-\alpha _{1}+\frac{n}{u}>0  \label{cond1}
	\end{equation}%
	and%
	\begin{equation}
	\sigma -\frac{n}{v}=-(1-\theta )\frac{n}{u}+\theta \Big(s-\frac{n}{p}\Big)%
	+\alpha _{1}-\big((1-\theta )\alpha _{2}+\theta \alpha _{3}\big),\quad 0<\theta <1.
	\label{cond1.2}
	\end{equation}%
	Assume that  $s>\sigma _{p,\beta }$ in the $\dot{K}F$%
	-case. $\newline
	\mathrm{(i)}$ Let $\alpha _{1}\leq \alpha _{2}\leq \alpha _{3}$. There is a constant $c>0$ such that for all $f\in \dot{K}_{u}^{\alpha _{2},\delta }\cap 
	\dot{K}_{p}^{\alpha _{3},\delta _{1}}B_{\beta }^{s},$%
	\begin{equation}
	\big\|f\big\|_{\dot{K}_{v}^{\alpha _{1},r}\dot{F}_{2}^{\sigma }}\leq c\big\|f%
	\big\|_{\dot{K}_{u}^{\alpha _{2},\delta }}^{1-\theta }\big\|f\big\|_{\dot{K}%
		_{p}^{\alpha _{3},\delta _{1}}\dot{B}_{\beta }^{s}}^{\theta }
	\label{key-est1}
	\end{equation}%
	with%
	\begin{equation*}
	\delta =\left\{ 
	\begin{array}{ccc}
	r, & \text{if} & \alpha _{2}=\alpha _{1}, \\ 
	\tau , & \text{if} & \alpha _{2}>\alpha _{1}.%
	\end{array}%
	\right. \quad \text{and}\quad \delta _{1}=\left\{ 
	\begin{array}{ccc}
	r, & \text{if} & \alpha _{3}=\alpha _{1}, \\ 
	\varrho , & \text{if} & \alpha _{3}>\alpha _{1}.%
	\end{array}%
	\right.
	\end{equation*}%
	$\mathrm{(ii)}$ Let $\frac{1}{r}\leq (1-\theta )\frac{n}{u}+\theta \frac{n}{p}$ and 
	\begin{equation*}
	\alpha _{1}=(1-\theta )\alpha _{2}+\theta \alpha _{3}.
	\end{equation*}%
	There is a constant $c>0$ such that for all $f\in \dot{K}_{u}^{\alpha
		_{2},u}F_{\infty }^{0}\cap \dot{K}_{p}^{\alpha _{3},p}A_{\infty }^{s},$ 
	\begin{equation*}
	\big\|f\big\|_{\dot{k}_{v,\sigma }^{\alpha _{1},r}}\leq c\big\|f\big\|_{\dot{%
			K}_{u}^{\alpha _{2},u}F_{\infty }^{0}}^{1-\theta }\big\|f\big\|_{\dot{K}%
		_{p}^{\alpha _{3},p}A_{\infty }^{s}}^{\theta }.
	\end{equation*}
\end{theorem}

\begin{proof}

	\textit{Proof of }(i)\textit{.} For technical reasons, we split
	the proof into two steps.
	
	\textit{Step 1.} We consider the case $p\leq u$. Let 
	\begin{equation*}
	f=\sum_{j=0}^{\infty }\mathcal{F}^{-1}\varphi _{j}\ast f,\quad f\in \mathcal{%
		S}^{\prime }(\mathbb{R}^{n}).
	\end{equation*}%
	Then it follows that%
	\begin{align*}
	f &=\sum_{j=0}^{J}\mathcal{F}^{-1}\varphi _{j}\ast f+\sum_{j=J+1}^{\infty }%
	\mathcal{F}^{-1}\varphi _{j}\ast f \\
	&=Q_{J}f+\sum_{j=J+1}^{\infty }\mathcal{F}^{-1}\varphi _{j}\ast f,\quad
	J\in \mathbb{N}.
	\end{align*}%
	Hence%
	\begin{equation}
	\big\|f\big\|_{\dot{k}_{v,\sigma }^{\alpha _{1},r}}\leq \big\|Q_{J}f\big\|_{%
		\dot{k}_{v,\sigma }^{\alpha _{1},r}}+\Big\|\sum_{j=J+1}^{\infty }\mathcal{F}%
	^{-1}\varphi _{j}\ast f\Big\|_{\dot{k}_{v,\sigma }^{\alpha _{1},r}}.
	\label{est-f}
	\end{equation}%
	Using Proposition \ref{Triebel1}, it follows that%
	\begin{equation}
	\big\|Q_{J}f\big\|_{\dot{k}_{v,\sigma }^{\alpha _{1},r}}\lesssim 2^{J(\frac{n%
		}{u}-\frac{n}{v}+\alpha _{2}-\alpha _{1}+\sigma )}\big\|f\big\|_{\dot{K}%
		_{u}^{\alpha _{2},\delta }}.\label{result1}
	\end{equation}%
	From the embedding%
	\begin{equation}
	\dot{K}_{v}^{\alpha _{1},r}B_{1}^{\sigma }\hookrightarrow \dot{k}_{v,\sigma
	}^{\alpha _{1},r},  \label{embedding}
	\end{equation}%
	see \eqref{aux7}, the last norm in \eqref{est-f} can be estimated by%
	\begin{align}
	c	\sum_{j=J+1}^{\infty }2^{j\sigma }\big\|\mathcal{F}^{-1}\varphi _{j}\ast f%
	\big\|_{\dot{K}_{v}^{\alpha _{1},r}} &\lesssim \sum_{j=J+1}^{\infty }2^{j(%
		\frac{n}{p}-\frac{n}{v}+\alpha _{3}-\alpha _{1}+\sigma )}\big\|\mathcal{F}%
	^{-1}\varphi _{j}\ast f\big\|_{\dot{K}_{p}^{\alpha _{3},\delta _{1}}} \nonumber \\
	&\lesssim 2^{J(\frac{n}{p}-\frac{n}{v}+\alpha _{3}-\alpha _{1}-s+\sigma )}%
	\big\|f\big\|_{\dot{K}_{p}^{\alpha _{3},\delta _{1}}B_{\beta }^{s}},\label{result2}
	\end{align}
	by Lemma \ref{Bernstein-Herz-ine1}, where the last estimate follows by %
	\eqref{cond1}. Plug \eqref{result1} and \eqref{result2} into \eqref{est-f} we obtain
	\begin{align*}
	\big\|f\big\|_{\dot{k}_{v,\sigma }^{\alpha _{1},r}} &\lesssim 2^{J(\frac{n}{%
			u}-\frac{n}{v}+\alpha _{2}-\alpha _{1}+\sigma )}\big\|f\big\|_{\dot{K}%
		_{u}^{\alpha _{2},\delta }}+2^{J(\frac{n}{p}-\frac{n}{v}+\alpha _{3}-\alpha
		_{1}-s+\sigma )}\big\|f\big\|_{\dot{K}_{p}^{\alpha _{3},\delta _{1}}B_{\beta
		}^{s}} \\
	&=c2^{J(\frac{n}{u}-\frac{n}{v}+\alpha _{2}-\alpha _{1}+\sigma )}\left( %
	\big\|f\big\|_{\dot{K}_{u}^{\alpha _{2},\delta }}+2^{J(\frac{n}{p}-\frac{n}{u%
		}-s-\alpha _{2}+\alpha _{3})}\big\|f\big\|_{\dot{K}_{p}^{\alpha _{3},\delta
			_{1}}B_{\beta }^{s}}\right) ,
	\end{align*}%
	with some positive constant $c$ independent of $J$. Again from, Lemma \ref%
	{Bernstein-Herz-ine1}, it follows that 
	\begin{equation}
	\dot{K}_{p}^{\alpha _{3},\delta _{1}}B_{\beta }^{s}\hookrightarrow \dot{K}%
	_{u}^{\alpha _{2},\delta }\label{key-emb},
	\end{equation}%
	since $s-\frac{n}{p}+\frac{n}{u}+\alpha _{2}-\alpha _{3}>0$. We choose $J\in 
	\mathbb{N}$ such that 
	\begin{equation*}
	2^{J(\frac{n}{p}-\frac{n}{u}-s-\alpha _{2}+\alpha _{3})}\approx \big\|f\big\|%
	_{\dot{K}_{u}^{\alpha _{2},\delta }}\big\|f\big\|_{\dot{K}_{p}^{\alpha
			_{3},\delta _{1}}B_{\beta }^{s}}^{-1}.
	\end{equation*}%
	We obtain\ 
	\begin{equation*}
	\big\|f\big\|_{\dot{k}_{v,\sigma }^{\alpha _{1},r}}\lesssim \big\|f\big\|_{%
		\dot{K}_{u}^{\alpha _{2},\delta }}^{1-\theta }\big\|f\big\|_{\dot{K}%
		_{p}^{\alpha _{3},\delta _{1}}B_{\beta }^{s}}^{\theta }.
	\end{equation*}%
	By \eqref{cond1} one has $s>\max \big(\sigma _{p},\alpha _{3}-n+\frac{n}{p}%
	\big)$ and by the fact that $-\frac{n}{u}<\alpha _{2}<n-\frac{n}{u},$%
	\begin{equation*}
	\sigma >\max \Big(0,\alpha _{1}+\frac{n}{v}-n\Big)
	\end{equation*}%
	and Theorem \ref{means-diff-cha}, or Proposition \ref{Triebel1 copy(3)}, can
	be used. Therefore 
	\begin{equation*}
	\big\|f\big\|_{\dot{K}_{v}^{\alpha _{1},r}\dot{F}_{2}^{\sigma }}\lesssim \big\|f\big\|_{\dot{k}_{v,\sigma }^{\alpha _{1},r}}
	\end{equation*}%
	and%
	\begin{equation*}
	\big\|f\big\|_{\dot{K}_{v}^{\alpha _{1},r}\dot{F}_{2}^{\sigma }}\lesssim %
	\big\|f\big\|_{\dot{K}_{u}^{\alpha _{2},\delta }}^{1-\theta }\left( \big\|f%
	\big\|_{\dot{K}_{p}^{\alpha _{3},\delta _{1}}}+\big\|f\big\|_{\dot{K}%
		_{p}^{\alpha _{3},\delta _{1}}\dot{B}_{\beta }^{s}}\right) ^{\theta }.
	\end{equation*}%
	In this estimate replace $f$ by $f(\lambda \cdot )$ we obtain%
	\begin{equation*}
	\big\|f\big\|_{\dot{K}_{v}^{\alpha _{1},r}\dot{F}_{2}^{\sigma }}\lesssim %
	\big\|f\big\|_{\dot{K}_{u}^{\alpha _{2},\delta }}^{1-\theta }\left( \lambda
	^{- s}\big\|f\big\|_{\dot{K}_{p}^{\alpha _{3},\delta _{1}}}+\big\|f%
	\big\|_{\dot{K}_{p}^{\alpha _{3},\delta _{1}}\dot{B}_{\beta }^{s}}\right)
	^{\theta }.
	\end{equation*}%
	Taking $\lambda $ large enough we obtain \eqref{key-est1} but with $p\leq u$.
	
	\textit{Step 2}. We consider the case $u<p$.  Taking $\lambda >0$ large enough such that%
	\begin{equation}
	\frac{\big\|f(\lambda \cdot )\big\|_{\dot{K}_{u}^{\alpha _{2},\delta }}}{%
		\big\|f(\lambda \cdot )\big\|_{\dot{K}_{p}^{\alpha _{3},\delta _{1}}B_{\beta
			}^{s}}}\leq 1 \label{key-emb1},
	\end{equation}%
	which is possible because of $s-\frac{n}{p}+\frac{n}{u}+\alpha _{2}-\alpha
	_{3}>0$, see \eqref{delation1}. As in Step 1, with  $f(\lambda \cdot )$ in place of $f$ and \eqref{key-emb1} in place of \eqref{key-emb}, we obtain the desired estimate. The proof of (i) is complete.
	
	\textit{Proof of }(ii)\textit{. }Observe that%
	\begin{equation*}
	\frac{n}{v_{1}}=\frac{n}{v}+\theta s-\sigma =(1-\theta )\frac{n}{u}+\theta 
	\frac{n}{p}
	\end{equation*}%
	and $\frac{\sigma }{s}\leq \theta <1$. Therefore 
	\begin{equation*}
	\dot{K}_{v_{1}}^{\alpha _{1},r}F_{\infty }^{\theta s}\hookrightarrow \dot{k}%
	_{v,\sigma }^{\alpha _{1},r},
	\end{equation*}%
	see Theorems \ref{embeddings3}. From \eqref{coincidence1}, %
	\eqref{coincidence3} and \eqref{Interpolation}, we obtain 
	\begin{equation*}
	\big\|f\big\|_{\dot{K}_{v_{1}}^{\alpha _{1},r}F_{\infty }^{\theta s}}\leq %
	\big\|f\big\|_{\dot{K}_{u}^{\alpha _{2},u}F_{\infty }^{0}}^{1-\theta }\big\|f%
	\big\|_{\dot{K}_{p}^{\alpha _{3},p}F_{\infty }^{\theta s}}^{\theta }.
	\end{equation*}%
	We have%
	\begin{equation*}
	\dot{K}_{p}^{\alpha _{3},p}A_{\beta }^{s}\hookrightarrow \dot{K}_{p}^{\alpha
		_{3},p}F_{\infty }^{\theta s}.
	\end{equation*}%
	This finishes the proof of (ii). The proof is complete. 
\end{proof}

\begin{remark}
	{\rm	$\mathrm{(i)}$ Taking $\alpha _{1}=\alpha _{2}=\alpha _{3}=0$ and $r=v$\ we
		obtain 
		\begin{align*}
		\big\|f\big\|_{\dot{H}_{v}^{\sigma }} &\leq c\big\|f\big\|_{\dot{K}%
			_{u}^{0,v}}^{1-\theta }\big\|f\big\|_{\dot{K}_{p}^{0,v}\dot{B}_{\beta
			}^{s}}^{\theta } \\
		&\leq c\big\|f\big\|_{u}^{1-\theta }\big\|f\big\|_{\dot{B}_{p,\beta
			}^{s}}^{\theta }
		\end{align*}%
		for all $f\in L_{u}\cap B_{p,\beta }^{s}$, because of $L_{u}\hookrightarrow 
		\dot{K}_{u}^{0,v}$ and $\dot{B}_{p,\beta }^{s}=\dot{K}_{p}^{0,p}\dot{B}%
		_{p,\beta }^{s}\hookrightarrow \dot{K}_{p}^{0,v}\dot{B}_{\beta }^{s}$, which
		has been proved by Triebel in \cite[Theorem 4.6]{Triebel13}.$\newline
		\mathrm{(ii)}$ Under the hypothesis of Theorem \ref{Triebel2}/(ii), with $%
		0<p<\frac{n}{s-\frac{\sigma }{\theta }}$ and $\frac{1}{r}\leq (1-\theta )%
		\frac{n}{u}+\theta (\frac{n}{p}-s+\frac{\sigma }{\theta })$, we have 
		\begin{equation*}
		\big\|f\big\|_{\dot{k}_{v,\sigma }^{\alpha _{1},r}}\leq c\big\|f\big\|_{\dot{%
				K}_{u}^{\alpha _{2},u}F_{2}^{0}}^{1-\theta }\big\|f\big\|_{\dot{K}_{p}^{\alpha _{3},\frac{1}{\frac{n}{p}-s+\frac{\sigma }{\theta }}}A_{\kappa}^{s}}^{\theta }
		\end{equation*}%
		for all $f\in \dot{K}_{u}^{\alpha _{2},u}F_{2}^{0}\cap \dot{K}_{p}^{\alpha
			_{3},\frac{1}{\frac{n}{p}-s+\frac{\sigma }{\theta }}}A_{\kappa }^{s}$, where%
		\begin{equation*}
		\kappa =\left\{ 
		\begin{array}{ccc}
		\frac{1}{\frac{n}{p}-s+\frac{\sigma }{\theta }}, & \text{if} & A=B, \\ 
		\infty , & \text{if} & A=F.%
		\end{array}%
		\right.
		\end{equation*}%
		Indeed, observe that%
		\begin{equation*}
		\frac{n}{v}=(1-\theta )\frac{n}{u}+\theta \Big(\frac{n}{p}-s+\frac{\sigma }{%
			\theta }\Big)=(1-\theta )\frac{n}{u}+\theta \frac{n}{u_{1}}
		\end{equation*}%
		and $\frac{\sigma }{\theta }-s\leq 0$. Therefore, from \eqref{coincidence1}, %
		\eqref{coincidence3} and \eqref{Interpolation}, we obtain 
		\begin{equation*}
		\big\|f\big\|_{\dot{k}_{v,\sigma }^{\alpha _{1},r}}\leq \big\|f\big\|_{\dot{K%
			}_{u}^{\alpha _{2},u}F_{2}^{0}}^{1-\theta }\big\|f\big\|_{\dot{K}%
			_{u_{1}}^{\alpha _{3},\frac{1}{\frac{n}{p}-s+\frac{\sigma }{\theta }}}F_{2}^{%
				\frac{\sigma }{\theta }}}^{\theta }.
		\end{equation*}%
		The result follows by 
		\begin{equation*}
		\dot{K}_{p}^{\alpha _{3},\frac{1}{\frac{n}{p}-s+\frac{\sigma }{\theta }}%
		}A_{\kappa }^{s}\hookrightarrow \dot{K}_{u_{1}}^{\alpha _{3},\frac{1}{\frac{n%
				}{p}-s+\frac{\sigma }{\theta }}}F_{2}^{\frac{\sigma }{\theta }},
		\end{equation*}%
		see Theorems \ref{embeddings3} and \ref{Franke-emb}.}
\end{remark}

\begin{theorem}
	\label{Triebel2.1}Let $\alpha _{1},\alpha _{2},\alpha _{3}\in \mathbb{R}%
	,0<p,\tau ,\beta ,\varrho \leq \infty $, $1<r,v,u<\infty ,$%
	\begin{equation*}
	s-\frac{n}{p}+\frac{n}{u}+\alpha _{2}-\alpha _{3}>-\frac{n}{v}+\alpha
	_{2}-\alpha _{1}+\frac{n}{u}>0
	\end{equation*}%
	and%
	\begin{equation*}
	\frac{n}{v}=(1-\theta )\frac{n}{u}+\theta \Big(\frac{n}{p}-s\Big)-\alpha
	_{1}+(1-\theta )\alpha _{2}+\theta \alpha _{3},\quad 0<\theta <1.
	\end{equation*}%
	Assume that $0<p,\tau <\infty $ and $s>\sigma _{p,\beta }$ in the $\dot{K}F$%
	-case. $\newline
	$Let $\delta $ and $\delta _{1}$ be as in Theorem \ref{Triebel2}/(i). Let $%
	\alpha _{1}\leq \alpha _{2}\leq \alpha _{3},v\geq \max (u,p),\alpha _{1}>-%
	\frac{n}{v},-\frac{n}{u}<\alpha _{2}<n-\frac{n}{u}$\ and $\alpha _{3}>-\frac{n}{p}$.
	We have 
	\begin{equation*}
	\big\|f\big\|_{\dot{K}_{v}^{\alpha _{1},r}}\lesssim \big\|f\big\|_{\dot{K}%
		_{u}^{\alpha _{2},\delta }}^{1-\theta }\big\|f\big\|_{\dot{K}_{p}^{\alpha
			_{3},\delta _{1}}A_{\beta }^{s}}^{\theta },  \label{Niremberg}
	\end{equation*}%
	holds for all $f\in \dot{K}_{u}^{\alpha _{2},\delta }\cap \dot{K}%
	_{p}^{\alpha _{3},\delta _{1}}A_{\beta }^{s}.$
\end{theorem}

\begin{proof}
	We employ the same notation and conventions as in Theorem %
	\ref{Triebel2}. As in Proposition \ref{Triebel1} 
	\begin{equation*}
	\big\|Q_{J}f\big\|_{\dot{K}_{v}^{\alpha _{1},r}}\lesssim 2^{J(\frac{n}{u}-%
		\frac{n}{v}+\alpha _{2}-\alpha _{1})}\big\|f\big\|_{\dot{K}_{u}^{\alpha
			_{2},\delta }},\quad J\in \mathbb{N}.
	\end{equation*}%
	Therefore, 
	\begin{equation*}
	\big\|f\big\|_{\dot{K}_{v}^{\alpha _{1},r}}\lesssim 2^{J(\frac{n}{u}-\frac{n%
		}{v}+\alpha _{2}-\alpha _{1})}\big\|f\big\|_{\dot{K}_{u}^{\alpha _{2},\delta
	}}+\sum_{j=J+1}^{\infty }\big\|\mathcal{F}^{-1}\varphi _{j}\ast f\big\|_{%
		\dot{K}_{v}^{\alpha _{1},r}},\quad J\in \mathbb{N}.
	\end{equation*}%
	Repeating the same arguments of Theorem \ref{Triebel2} we obtain the desired
	estimate.
\end{proof}

\begin{remark}
	{\rm	Under the same hypothesis of Theorem \ref{Triebel2.1}, with $1<p<\infty ,-%
		\frac{n}{p}<\alpha _{3}<n-\frac{n}{p},r=v\ $and $\beta =2$, we obtain%
		\begin{align*}
		\big\||\cdot |^{\alpha _{1}}f\big\|_{v} &\lesssim \big\|f\big\|_{\dot{K}%
			_{u}^{\alpha _{2},v}}^{1-\theta }\big\|f\big\|_{\dot{K}_{p}^{\alpha
				_{3},v}F_{2}^{s}}^{\theta } \\
		&\lesssim \big\||\cdot |^{\alpha _{2}}f\big\|_{u}^{1-\theta }\big\|f\big\|_{%
			\dot{k}_{p,s}^{\alpha _{3},v}}^{\theta } \\
		&\lesssim \big\||\cdot |^{\alpha _{2}}f\big\|_{u}^{1-\theta }\big\|f\big\|_{%
			\dot{k}_{p,s}^{\alpha _{3},p}}^{\theta }
		\end{align*}%
		for all $f\in L^{u}(\mathbb{R}^{n},|\cdot |^{\alpha _{2}u})\cap \dot{k}_{p,s}^{\alpha
			_{3},p} $, because of 
		\begin{equation*}
		\dot{K}_{u}^{\alpha _{2},u}\hookrightarrow \dot{K}_{u}^{\alpha _{2},v}\quad 
		\text{and}\quad \dot{k}_{p,s}^{\alpha _{3},p}\hookrightarrow \dot{k}%
		_{p,s}^{\alpha _{3},v}.
		\end{equation*}%
		In particular if $s=m\in \mathbb{N}$, then we obtain%
		\begin{align}
		\big\||\cdot |^{\alpha _{1}}f\big\|_{v} &\lesssim \big\|f\big\|_{\dot{K}%
			_{u}^{\alpha _{2},v}}^{1-\theta }\Big(\sum\limits_{|\beta |\leq m}\Big\|%
		\frac{\partial ^{\beta }f}{\partial ^{\beta }x}\Big\|_{\dot{K}_{p}^{\alpha
				_{3},v}}\Big)^{\theta }  \label{conclutions} \\
		&\lesssim \big\||\cdot |^{\alpha _{2}}f\big\|_{u}^{1-\theta }\Big(%
		\sum\limits_{|\beta |\leq m}\Big\||\cdot |^{\alpha _{3}}\frac{\partial
			^{\beta }f}{\partial ^{\beta }x}\Big\|_{p}\Big)^{\theta }  \notag
		\end{align}%
		for all $f\in L^{u}(\mathbb{R}^{n},|\cdot |^{\alpha _{2}u})\cap W_{p}^{m}(\mathbb{R}^{n},|\cdot |^{\alpha_{3}u})$. As in \cite[Theorem 4.6]{Triebel13} replace $f$ in %
		\eqref{conclutions} by $f(\lambda \cdot )$ with $\lambda >0$, the sum $%
		\sum\limits_{|\beta |\leq m}\cdot \cdot \cdot $ can be replaced by$%
		\sum\limits_{0<|\beta |\leq m}\cdot \cdot \cdot .$}
\end{remark}

From Proposition \ref{Triebel1 copy(2)} and Theorem \ref{Triebel2}/(i) we obtain the following statement.

\begin{theorem}
	\label{Triebel2 copy(1)}Let $1<p,\varrho <\infty ,0<\tau \leq \infty $, $%
	1<r,v,u<\infty ,\sigma \geq 0,$ \eqref{cond1.1}, \eqref{cond1} and %
	\eqref{cond1.2} with $\alpha _{3}<n-\frac{n}{p}$. Let $\alpha _{1}\leq
	\alpha _{2}\leq \alpha _{3}$. There is a constant $c>0$ such that for all $%
	f\in \dot{K}_{u}^{\alpha _{2},\delta }\cap \dot{k}_{p,s}^{\alpha _{3},\delta
		_{1}},$%
	\begin{equation*}
	\big\|(-\Delta ^{\frac{\sigma }{2}})f\big\|_{\dot{K}_{v}^{\alpha
			_{1},r}}\leq c\big\|f\big\|_{\dot{K}_{u}^{\alpha _{2},\delta }}^{1-\theta }%
	\big\|(-\Delta ^{\frac{s}{2}})f\big\|_{\dot{K}_{p}^{\alpha _{3},\delta
			_{1}}}^{\theta }
	\end{equation*}%
	with%
	\begin{equation*}
	\delta =\left\{ 
	\begin{array}{ccc}
	r, & \text{if} & \alpha _{2}=\alpha _{1}, \\ 
	\tau , & \text{if} & \alpha _{2}>\alpha _{1}.%
	\end{array}%
	\right. \quad \text{and}\quad \delta _{1}=\left\{ 
	\begin{array}{ccc}
	r, & \text{if} & \alpha _{3}=\alpha _{1}, \\ 
	\varrho , & \text{if} & \alpha _{3}>\alpha _{1}.%
	\end{array}%
	\right.
	\end{equation*}
\end{theorem}

In the next we study the case when $p\leq v< u$ in
Theorem \ref{Triebel2}.

\begin{theorem}
	\label{Triebel3}Let $0<p,\tau <\infty ,0<\beta ,\kappa \leq \infty
	,1<r,v<\infty ,\sigma \geq 0,1< u<\infty ,$%
	\begin{equation*}
	-\frac{n}{v}<\alpha _{1}<n-\frac{n}{v},\quad -\frac{n}{u}<\alpha _{2}<n-\frac{n}{u},\quad
	\alpha _{3}>-\frac{n}{p},
	\end{equation*}%
	\begin{equation*}
	s-\frac{n}{p}+\frac{n}{u}+\alpha _{2}-\alpha _{3}>\sigma -\frac{n}{v}+\alpha
	_{2}-\alpha _{1}+\frac{n}{u}>0
	\end{equation*}%
	and%
	\begin{equation*}
	\sigma -\frac{n}{v}=-(1-\theta )\frac{n}{u}+\theta \Big(s-\frac{n}{p}\Big)%
	+\alpha _{1}-((1-\theta )\alpha _{2}+\theta \alpha _{3}),\quad 0<\theta <1.
	\end{equation*}$%
	\newline
	\mathrm{(i)}$ Let $p\leq v<u,\alpha _{2}-\alpha _{1}>\frac{n}{v}-\frac{n}{u}%
	\ $and$\ \alpha _{3}=\alpha _{2}$. There is a constant $c>0$ such that for
	all $f\in \dot{K}_{u}^{\alpha _{2},\tau }\cap \dot{K}_{p}^{\alpha _{3},\tau
	}F_{\beta }^{s},$%
	\begin{equation}
	\big\|f\big\|_{\dot{k}_{v,\sigma }^{\alpha _{1},r}}\leq c\big\|f\big\|_{\dot{%
			K}_{u}^{\alpha _{2},\tau }}^{1-\theta }\big\|f\big\|_{\dot{K}_{p}^{\alpha
			_{3},\tau }F_{\beta }^{s}}^{\theta }.  \label{CKN1}
	\end{equation}%
	$\mathrm{(ii)}$ Let $p\leq v<u,\alpha _{2}-\alpha _{1}>\frac{n}{v}-\frac{n}{u%
	}\ $and$\ \alpha _{3}>\alpha _{2}$. There is a constant $c>0$ such that %
	\eqref{CKN1} holds for all $f\in \dot{K}_{u}^{\alpha _{2},\tau }\cap \dot{K}%
	_{p}^{\alpha _{3},\kappa }F_{\beta }^{s}$ with $\dot{K}_{p}^{\alpha
		_{3},\kappa }F_{\beta }^{s}$ in place of $\dot{K}_{p}^{\alpha _{3},\tau
	}F_{\beta }^{s}$.
	
\end{theorem}

\begin{proof}
	Recall that, as in Theorem \ref{Triebel2}, one has the estimate 
	\begin{equation*}
	\big\|f\big\|_{\dot{k}_{v,\sigma }^{\alpha _{1},r}}\leq \big\|Q_{J}f\big\|_{%
		\dot{k}_{v,\sigma }^{\alpha _{1},r}}+\Big\|\sum_{j=J+1}^{\infty }\mathcal{F}%
	^{-1}\varphi _{j}\ast f\Big\|_{\dot{k}_{v,\sigma }^{\alpha _{1},r}},\quad
	J\in \mathbb{N}.
	\end{equation*}%
	From Proposition \ref{Triebel1}/(ii),%
	\begin{equation*}
	\big\|Q_{J}f\big\|_{\dot{k}_{v,\sigma }^{\alpha _{1},r}}\leq c2^{J(\frac{n}{u%
		}-\frac{n}{v}+\alpha _{2}-\alpha _{1}+\sigma )}\big\|f\big\|_{\dot{K}%
		_{u}^{\alpha _{2},\tau }},
	\end{equation*}%
	which is possible since 
	\begin{equation*}
	\frac{n}{v}+\alpha _{1}-\alpha _{2}\leq \frac{n}{u}< \frac{n}{v}.
	\end{equation*}%
	Using again the embedding \eqref{embedding} and Lemma \ref%
	{Bernstein-Herz-ine1}, we get%
	\begin{align*}
	\Big\|\sum_{j=J+1}^{\infty }\mathcal{F}^{-1}\varphi _{j}\ast f\Big\|_{\dot{k}%
		_{v,\sigma }^{\alpha _{1},r}} &\lesssim \sum_{j=J+1}^{\infty }2^{j\sigma }%
	\big\|\mathcal{F}^{-1}\varphi _{j}\ast f\big\|_{\dot{K}_{v}^{\alpha _{1},r}}
	\\
	&\lesssim \sum_{j=J+1}^{\infty }2^{j(\frac{n}{p}-\frac{n}{v}+\alpha
		_{3}-\alpha _{1}+\sigma )}\big\|\mathcal{F}^{-1}\varphi _{j}\ast f\big\|_{%
		\dot{K}_{p}^{\alpha _{3},\vartheta }},
	\end{align*}%
	where%
	\begin{equation*}
	\vartheta =\left\{ 
	\begin{array}{ccc}
	\tau , & \text{if} & \alpha _{3}=\alpha _{2}, \\ 
	\kappa , & \text{if} & \alpha _{3}>\alpha _{2}.%
	\end{array}%
	\right.
	\end{equation*}%
	Therefore, $\big\|f\big\|_{\dot{k}_{v,\sigma }^{\alpha _{1},r}}$ can be
	estimated by%
	\begin{align*}
	&c2^{J(\frac{n}{u}-\frac{n}{v}+\alpha _{2}-\alpha _{1}+\sigma )}\big\|f%
	\big\|_{\dot{K}_{u}^{\alpha _{2},\tau }}+2^{J(\frac{n}{p}-\frac{n}{v}+\alpha
		_{3}-\alpha _{1}-s+\sigma )}\big\|f\big\|_{\dot{K}_{p}^{\alpha
			_{3},\vartheta }F_{\beta }^{s}} \\
	&=c2^{J(\frac{n}{u}-\frac{n}{v}+\alpha _{2}-\alpha _{1}+\sigma )}\left( %
	\big\|f\big\|_{\dot{K}_{u}^{\alpha _{2},\tau }}+2^{J(\frac{n}{p}-\frac{n}{u}%
		-s-\alpha _{2}+\alpha _{3})}\big\|f\big\|_{\dot{K}_{p}^{\alpha
			_{3},\vartheta }F_{\beta }^{s}}\right) ,
	\end{align*}%
	where the positive constant $c>0$ is independent of $J$. Observe that%
	\begin{equation*}
	\dot{K}_{p}^{\alpha _{3},\vartheta }F_{\beta }^{s}\hookrightarrow \dot{K}%
	_{u}^{\alpha _{2},\tau },
	\end{equation*}%
	since $s-\frac{n}{p}+\frac{n}{u}+\alpha _{2}-\alpha _{3}>0$. We choose $J\in 
	\mathbb{N}$ such that 
	\begin{equation*}
	2^{J(\frac{n}{p}-\frac{n}{u}-s-\alpha _{2}+\alpha _{3})}\approx \big\|f\big\|%
	_{\dot{K}_{u}^{\alpha _{2},\tau }}\big\|f\big\|_{\dot{K}_{p}^{\alpha
			_{3},\vartheta }F_{\beta }^{s}}^{-1},
	\end{equation*}%
	we obtain the desired estimate. The proof is complete.
\end{proof}

As in Theorem\ \ref{Triebel2.1} combined with Theorem \ref{Triebel3} we
obtain the following conclusion.

\begin{theorem}
	\label{Triebel4.1}Under the hypothesis of Theorem \ref{Triebel3} with\ $%
	\alpha _{1}>-\frac{n}{v}$ and $\sigma =0$, we have the estimates with $\dot{K%
	}_{v}^{\alpha _{1},r}$ replaced by $\dot{k}_{v,\sigma }^{\alpha _{1},r}$.
\end{theorem}

Finally we study the case of $v\leq \min (p,u)$.

\begin{theorem}
	\label{Triebel4}Let $1<r<\infty ,0<p,\beta ,\tau \leq \infty ,1<v\leq \min
	(p,u),\alpha _{2}-\alpha _{1}>\frac{n}{v}-\frac{n}{\max (p,u)},\alpha
	_{3}\geq \alpha _{2},\sigma \geq 0$,%
	\begin{equation*}
	-\frac{n}{v}<\alpha _{1}<n-\frac{n}{v},\quad -\frac{n}{u}<\alpha _{2}<n-\frac{n}{u},\quad
	\alpha _{3}>-\frac{n}{p}
	\end{equation*}%
	and 
	\begin{equation*}
	s-\frac{n}{p}+\frac{n}{u}+\alpha _{2}-\alpha _{3}>\sigma -\frac{n}{v}+\alpha
	_{2}-\alpha _{1}+\frac{n}{u}>0.  \label{case1}
	\end{equation*}%
	Assume that $0<p,\tau <\infty $ and $s>\sigma _{p,\beta }$ in the $\dot{K}F$%
	-case. There is a constant $c>0$ such that for all $f\in \dot{K}_{u}^{\alpha
		_{2},\tau }\cap \dot{K}_{p}^{\alpha _{3},\tau }A_{\beta }^{s},$%
	\begin{equation*}
	\big\|f\big\|_{\dot{k}_{v,\sigma }^{\alpha _{1},r}}\leq c\big\|f\big\|_{\dot{%
			K}_{u}^{\alpha _{2},\tau }}^{1-\theta }\big\|f\big\|_{\dot{K}_{p}^{\alpha
			_{3},\tau }A_{\beta }^{s}}^{\theta }
	\end{equation*}%
	with%
	\begin{equation*}
	\sigma -\frac{n}{v}=-(1-\theta )\frac{n}{u}+\theta \Big(s-\frac{n}{p}\Big)%
	+\alpha _{1}-\big((1-\theta )\alpha _{2}+\theta \alpha _{3}\big).
	\end{equation*}
\end{theorem}

\begin{proof}
	By similarity, we only consider $\dot{K}_{p}^{\alpha
		_{3},\tau }B_{\beta }^{s}$. We split the proof into two steps.
	
	\textit{Step 1}. We consider the case $p\leq u$. We employ the same notation as in Theorem \ref%
	{Triebel2}. In view of Theorem \ref{Triebel3} we need only to estimate%
	\begin{equation*}
	\Big\|\sum_{j=J+1}^{\infty }\mathcal{F}^{-1}\varphi _{j}\ast f\Big\|_{\dot{k}%
		_{v,\sigma }^{\alpha _{1},r}},\quad J\in \mathbb{N}.
	\end{equation*}%
	Using the embedding \eqref{embedding} and Lemma \ref{Bernstein-Herz-ine2},
	we obtain%
	\begin{align*}
	\Big\|\sum_{j=J+1}^{\infty }\mathcal{F}^{-1}\varphi _{j}\ast f\Big\|_{\dot{k}%
		_{v,\sigma }^{\alpha _{1},r}} &\lesssim \sum_{j=J+1}^{\infty }2^{j\sigma }%
	\big\|\mathcal{F}^{-1}\varphi _{j}\ast f\big\|_{\dot{K}_{v}^{\alpha _{1},r}}
	\\
	&\lesssim \sum_{j=J+1}^{\infty }2^{j(\frac{n}{p}-\frac{n}{v}+\alpha
		_{2}-\alpha _{1}+\sigma )}\big\|\mathcal{F}^{-1}\varphi _{j}\ast f\big\|_{%
		\dot{K}_{p}^{\alpha _{3},\tau }}.
	\end{align*}%
	with is possible since 
	\begin{equation*}
	\frac{n}{v}+\alpha _{1}-\alpha _{2}<\frac{n}{p}\leq \frac{n}{v}.
	\end{equation*}%
	Repeating the same arguments of Theorem \ref{Triebel2} we obtain the desired
	estimate.
	
	\textit{Step 2.} We consider the case $u<p$. Using a combination of the arguments used
	in the corresponding step of the proof of Theorem \ref{Triebel2} and those
	used in the first step above, we arrive at the desired estimate. 
\end{proof}

Similarly we obtain the following conclusion.

\begin{theorem}
	\label{Triebel5}Under the hypothesis of Theorem \ref{Triebel4} with\ $\sigma=0$, we have%
	\begin{equation*}
	\big\|f\big\|_{\dot{K}_{v}^{\alpha _{1},r}}\lesssim \big\|f\big\|_{\dot{K}%
		_{u}^{\alpha _{2},\tau }}^{1-\theta }\big\|f\big\|_{\dot{K}_{p}^{\alpha
			_{3},\tau }A_{\varrho }^{s}}^{\theta }
	\end{equation*}%
	for all $f\in \dot{K}_{u}^{\alpha _{2},\tau }\cap \dot{K}_{p}^{\alpha
		_{2},\tau }A_{\varrho }^{s}.$
\end{theorem}

\begin{remark}
	{\rm	Under the same hypothesis of Theorems \ref{Triebel4.1}\ and \ref{Triebel5},
		with\ $r=v,\sigma =0,\tau =\max (u,p)$ and $\beta =2$, we improve
		Caffarelli-Kohn-Nirenberg inequality \eqref{CKN} in some sense.}
\end{remark}

\subsection{ \large  CKN inequalities in Besov-Morrey and Triebel-Lizorkin-Morrey spaces}

In this section, we investigate the Caffarelli, Kohn and Nirenberg
inequalities in $\mathcal{E}_{p,q,u}^{s}$ and $\mathcal{N}_{p,q,u}^{s}$
spaces. The main results of this section based on the following Lemma.

\begin{lemma}
	\label{Polya-Nikolskij1}Let $1<u\leq p<\infty ,1<s\leq q<\infty $ and $%
	R>0$.\newline
	$\mathrm{(i)}$ \textit{Assume that }$1\leq v\leq u$. \textit{There exists a
		constant }$c>0$\textit{\ independent of }$R$\textit{\ such that for all }$%
	f\in M_{v}^{\frac{v}{u}p}\cap M_{s}^{q}$\textit{\ with}
	$\mathrm{supp}$\textit{\ }$\mathcal{F}f\subset
	\{\xi :|\xi |\leq R\}$, we have%
	\begin{equation*}
	\big\|f\big\|_{M_{u}^{p}}\leq cR^{\frac{n}{q}-\frac{vn}{qu}}\big\|f%
	\big\|_{M_{s}^{q}}^{1-\frac{v}{u}}\big\|f\big\|_{M_{v}^{%
			\frac{v}{u}p}}^{\frac{v}{u}}.  \label{PN1}
	\end{equation*}%
	$\mathrm{(ii)}$ \textit{Assume that }$\frac{u}{p}\leq \frac{s}{q}$ and $%
	q\leq p$. \textit{There exists a constant }$c>0$\textit{\ independent of }$R$%
	\textit{\ such that for all }$f\in M_{s}^{q}$\textit{\ with }$\mathrm{supp}$\textit{\ }$\mathcal{F}f\subset
	\{\xi :|\xi |\leq R\}$, we have%
	\begin{equation*}
	\big\|f\big\|_{M_{u}^{p}}\leq cR^{\frac{n}{q}-\frac{n}{p}}\big\|f%
	\big\|_{M_{s}^{q}}.  \label{PN2}
	\end{equation*}
\end{lemma}

\begin{proof}
	
	We split the proof in two steps.
	
	\textit{Step 1.} We will prove (i). Let $B$ be a ball of $\mathbb{R}^{n}$. Write 
	\begin{equation*}
	\big\|\left\vert B\right\vert ^{\frac{1}{p}-\frac{1}{u}}f\chi _{B}\big\|%
	_{u}^{u}=u\int_{0}^{\infty }t ^{u-1}|\{x\in B:|f(x)|\left\vert
	B\right\vert ^{\frac{1}{p}-\frac{1}{u}}>t \}|dt<\infty .
	\end{equation*}%
	We have 
	\begin{equation*}
	|f(x)|\leq cR^{\frac{n}{q}}\big\|f\big\|_{M_{s}^{q}},\quad x\in 
	\mathbb{R}^{n},
	\end{equation*}%
	see \cite[Proposition 2.1]{SST09} where $c>0$\ independent of $R$. Let $%
	p_{0}=\frac{v}{u}$. Clearly%
	\begin{align*}
	|f(x)| &=|f(x)|^{p_{0}}|f(x)|^{1-p_{0}} \\
	&\lesssim |f(x)|^{p_{0}}\big(R^{\frac{n}{q}}\big\|f\big\|_{M_{s}^{q}}\big)^{1-p_{0}} \\
	&=c|f(x)|^{p_{0}}d^{1-p_{0}},
	\end{align*}%
	which yields that%
	\begin{align*}
	\big\|\left\vert B\right\vert ^{\frac{1}{p}-\frac{1}{u}}f\chi _{B}\big\|%
	_{u}^{u} &\leq u\int_{0}^{\infty }t ^{u-1}|\{x\in B:|f(x)|\left\vert
	B\right\vert ^{\frac{1}{pp_{0}}-\frac{1}{v}}>cd^{1-\frac{1}{p_{0}}}t ^{%
		\frac{1}{p_{0}}}\}|dt  \\
	&=cud^{u-v}\int_{0}^{\infty }\lambda ^{v-1}|\{x\in B:|f(x)|\left\vert
	B\right\vert ^{\frac{u}{pv}-\frac{1}{v}}>\lambda \}|d\lambda  ,
	\end{align*}%
	after the change the variable $\lambda ^{p_{0}}c^{-p_{0}}d^{1-p_{0}}=t
	. $ The last expression is clearly bounded by%
	\begin{equation*}
	cd^{u-v}\big\|f\big\|_{M_{v}^{\frac{pv}{u}}}^{v}\leq cR^{n\frac{u-v%
		}{q}}\big\|f\big\|_{M_{v}^{\frac{pv}{u}}}^{v}\big\|f\big\|_{%
		M_{s}^{q}}^{u-v}.
	\end{equation*}
	
	\textit{Step 2.} We will prove (ii). If $p=q$, then $u\leq s$ and the estimate
	follows by the H\"{o}lder inequality. Assume that $q<p$ and we choose $v>0$ such that $\max(1,\frac{qu}{p})<v\leq u<\frac{pu}{q}$. By Step 1, we only need to estimate $R^{\frac{n}{q}-\frac{vn}{qu}}\big\|f\big\|_{M_{v}^{\frac{v}{u}p}}^{%
		\frac{v}{u}}$. Write%
	\begin{equation*}
	R^{\frac{n}{q}-\frac{vn}{qu}}\big\|f\big\|_{M_{v}^{\frac{v}{u}p}}^{%
		\frac{v}{u}}=R^{\frac{n}{q}-\frac{n}{p}}\big\|R^{\frac{nu}{pv}-\frac{n}{q}}f%
	\big\|_{M_{v}^{\frac{v}{u}p}}^{\frac{v}{u}}.
	\end{equation*}%
	Let $\{\varphi _{j}\}_{j\in \mathbb{N}_{0}}$\ be a
	resolution of unity. Observe that%
	\begin{equation*}
	\mathcal{F}^{-1}\varphi _{j}\ast f=0 \quad  {\rm if}\quad R<2^{j-1} ,\quad j\in \mathbb{N}_{0}.
	\end{equation*}%
	This observation together with \eqref{Morrey} yield
	\begin{align*}
	\big\|R^{\frac{nu}{pv}-\frac{n}{q}}f\big\|_{M_{v}^{\frac{v}{u}%
			p}} &\approx\Big\|\Big(\sum\limits_{j\in \mathbb{N}_{0},2^{j-1}\leq R}^{\infty }R^{\frac{2nu}{pv}-\frac{2n}{q}}\left\vert \mathcal{F}%
	^{-1}\varphi _{j}\ast f\right\vert ^{2}\Big)^{1/2}\Big\|_{M_{v}^{\frac{v}{u}p}}\\
	&\lesssim\big\|f\big\|_{\mathcal{E}_{\frac{v}{u}p,2,v}^{\frac{nu}{pv}-\frac{n}{q}}} \\
	&\lesssim \big\|f\big\|_{M_{s}^{q}},
	\end{align*}
	which follows by the Sobolev embedding, see Theorem \ref{Sobolev-embeddings},%
	\begin{equation*}
	M_{s}^{q}=\mathcal{E}_{q,2,s}^{0}\hookrightarrow \mathcal{E}_{%
		\frac{v}{u}p,2,v}^{\frac{nu}{pv}-\frac{n}{q}},
	\end{equation*}%
	since%
	\begin{equation*}
	-\frac{n}{q}=\frac{nu}{pv}-\frac{n}{q}-\frac{nu}{pv},\quad q<\frac{vp}{u}%
	\quad \text{and}\quad \frac{u}{p}\leq \frac{s}{q}.\text{ }
	\end{equation*}%
	The lemma is proved. 
\end{proof}

\begin{proposition}
	\label{Triebel1 copy(5)}\ Let\textit{\ }$1< u\leq p<\infty ,1<q<
	\infty $ and $s>0$.\newline
	$\mathrm{(i)}$ Let $f\in \mathcal{N}_{p,q,u}^{s}$.\textit{\ Then}%
	\begin{equation}
	\big\|f\big\|_{\mathcal{N}_{p,q,u}^{s}}\approx \big\|f\big\|_{M_{u}^{p}}+\big\|f\big\|_{\mathcal{\dot{N}}_{p,q,u}^{s}},  \label{new-norm3}
	\end{equation}%
	where 
	\begin{equation*}
	\big\|f\big\|_{\mathcal{\dot{N}}_{p,q,u}^{s}}=\Big\|\Big(\sum\limits_{j=-%
		\infty }^{\infty }2^{qjs}\left\vert \mathcal{F}^{-1}\varphi ^{j}\ast
	f\right\vert ^{q}\Big)^{1/q}\Big\|_{M_{u}^{p}}.
	\end{equation*}%
	$\mathrm{(ii)}$ Let\ $f\in \mathcal{E}_{p,q,u}^{s}$\textit{. Then}%
	\begin{equation}
	\big\|f\big\|_{\mathcal{E}_{p,q,u}^{s}}\approx \big\|f\big\|_{M_{u}^{p}}+\big\|f\big\|_{\mathcal{\dot{E}}_{p,q,u}^{s}}  \label{new-norm4}
	\end{equation}%
	where 
	\begin{equation*}
	\big\|f\big\|_{\mathcal{\dot{E}}_{p,q,u}^{s}}=\Big\|\Big(\sum\limits_{j=-%
		\infty }^{\infty }2^{qjs}\left\vert \mathcal{F}^{-1}\varphi ^{j}\ast
	f\right\vert ^{q}\Big)^{1/q}\Big\|_{M_{u}^{p}}.
	\end{equation*}
\end{proposition}
\begin{proof}
	By similarity, we prove only (ii).\textbf{\ }We have as in
	the proof of Proposition \ref{Triebel1 copy(3)} that\textbf{\ } 
	\begin{equation*}
	\big\|f\big\|_{\mathcal{\dot{E}}_{p,q,u}^{s}}\lesssim \big\|f\big\|_{\mathcal{E}_{p,q,u}^{s}}.
	\end{equation*}%
	The only difference with the proof of Proposition \ref{Triebel1 copy(3)}
	consists in the fact that we use \cite[Lemma 2.5]{TangXu05}. Since $s>0$ we
	observe%
	\begin{equation*}
	\big\|f\big\|_{M_{u}^{p}}\approx \big\|f\big\|_{\mathcal{E}%
		_{p,2,u}^{0}}\lesssim \big\|f\big\|_{\mathcal{E}_{p,q,u}^{s}}.
	\end{equation*}
	Now we prove the opposite inequality. Obviously $\big\|\mathcal{F}%
	^{-1}\varphi _{0}\ast f\big\|_{M_{u}^{p}}$ can be estimated from
	above by $\big\|f\big\|_{M_{u}^{p}}$, which completes the proof. . 
\end{proof}

\begin{theorem}
	\label{Triebel1.1}Let $1<u\leq p<\infty $ and $1<v\leq q<\infty $. \textit{%
		Assume that }$\frac{u}{p}\leq \frac{v}{q},q\leq p$\ and\ $\sigma \geq 0$. 
	\textit{Then for all }$f\in M_{v}^{q}$ and all $J\in \mathbb{N}$,%
	\begin{equation*}
	\big\|Q_{J}f\big\|_{\mathcal{E}_{p,2,u}^{\sigma }}\leq c2^{Jn(\frac{1}{q}-%
		\frac{1}{p})+\sigma }\big\|f\big\|_{M_{v}^{q}},
	\end{equation*}%
	\textit{where the positive constant }$c$\textit{\ is independent of }$J$%
	\textit{.}
\end{theorem}
\begin{proof}
	Let $\sigma =\theta m+(1-\theta )0$, $\alpha \in \mathbb{N}%
	^{n}$ with $0<\theta <1$ and $|\alpha |\leq m$. We have%
	\begin{equation*}
	\big\|Q_{J}f\big\|_{\mathcal{E}_{p,2,u}^{\sigma }}\leq \big\|Q_{J}f\big\|_{%
		\mathcal{E}_{p,2,u}^{0}}^{1-\theta }\big\|Q_{J}f\big\|_{\mathcal{E}%
		_{p,2,u}^{m}}^{\theta }.
	\end{equation*}%
	Observe that%
	\begin{equation*}
	\mathcal{E}_{p,2,u}^{m}=M_{u}^{m,p}\quad \text{and}\quad \mathcal{E%
	}_{p,2,u}^{0}=M_{u}^{p},
	\end{equation*}%
	which yield that%
	\begin{equation*}
	\big\|Q_{J}f\big\|_{\mathcal{E}_{p,2,u}^{\sigma }}\leq \big\|Q_{J}f\big\|_{%
		M_{u}^{p}}^{1-\theta }\big\|Q_{J}f\big\|_{M
		_{u}^{m,p}}^{\theta },
	\end{equation*}%
	where the positive constant $c$\ is independent of $J$. Lemma \ref%
	{Polya-Nikolskij1} yields that%
	\begin{equation*}
	\big\|D^{\alpha }(Q_{J}f)\big\|_{M_{u}^{p}}\lesssim 2^{Jn(\frac{1}{%
			q}-\frac{1}{p})+|\alpha |}\big\|f\big\|_{M_{v}^{q}}.
	\end{equation*}%
	Therefore,%
	\begin{equation*}
	\big\|Q_{J}f\big\|_{\mathcal{E}_{p,2,u}^{\sigma }}\lesssim 2^{Jn(\frac{1}{q}-%
		\frac{1}{p})+\sigma }\big\|f\big\|_{M_{v}^{q}}.
	\end{equation*}%
	This finish the proof.
\end{proof}

Now we are in position to state the main result of this section.

\begin{theorem}
	\label{Triebel1.2}Let $1<u\leq p<\infty ,1<\mu \leq \delta <\infty ,1<\beta
	<\infty $, $\sigma \geq 0$ and $1<v\leq q<\infty $. \textit{Assume that }%
	\begin{equation*}
	\frac{u}{p}\leq \frac{\mu }{\delta }\leq \frac{v}{q},\quad s>0\quad \text{and%
	}\quad p\geq \delta \geq q .
	\end{equation*}%
	Let 
	\begin{equation*}
	s-\frac{n}{q}>\sigma -\frac{n}{p}\quad \text{and}\quad \sigma -\frac{n}{p}%
	=-(1-\theta )\frac{n}{\delta }+\theta \Big(s-\frac{n}{q}\Big),\quad 0<\theta
	<1.  \label{condition1}
	\end{equation*}%
	\textit{Then }%
	\begin{equation}
	\big\|f\big\|_{\mathcal{\dot{E}}_{p,2,u}^{\sigma }}\lesssim \big\|f\big\|_{%
		M_{\mu }^{\delta }}^{1-\theta }\big\|f\big\|_{\mathcal{\dot{N}}%
		_{q,\beta ,v}^{s}}^{\theta },\mathcal{\quad }\sigma >0  \label{result-Morrey}
	\end{equation}%
	and%
	\begin{equation}
	\big\|f\big\|_{M_{u}^{p}}\lesssim \big\|f\big\|_{M_{\mu
		}^{\delta }}^{1-\theta }\big\|f\big\|_{\mathcal{\dot{N}}_{q,\beta
			,v}^{s}}^{\theta }  \label{result-Morrey1}
	\end{equation}%
	\textit{for all }$f\in M_{\mu }^{\delta }\cap \mathcal{N}_{q,\beta
		,v}^{s}$.
\end{theorem}

\begin{proof}
	We have%
	\begin{equation*}
	f=Q_{J}f+\sum_{j=J+1}^{\infty }\mathcal{F}^{-1}\varphi _{j}\ast f,\quad J\in 
	\mathbb{N}.
	\end{equation*}%
	Hence%
	\begin{equation}
	\big\|f\big\|_{\mathcal{E}_{p,2,u}^{\sigma }}\leq \big\|Q_{J}f\big\|_{%
		\mathcal{E}_{p,2,u}^{\sigma }}+\Big\|\sum_{j=J+1}^{\infty }\mathcal{F}%
	^{-1}\varphi _{j}\ast f\Big\|_{\mathcal{E}_{p,2,u}^{\sigma }}.
	\label{est-f1}
	\end{equation}%
	Using Theorem \ref{Triebel1.1}, it follows that%
	\begin{equation*}
	\big\|Q_{J}f\big\|_{\mathcal{E}_{p,2,u}^{\sigma }}\lesssim 2^{Jn(\frac{1}{%
			\delta }-\frac{1}{p})+\sigma J}\big\|f\big\|_{M_{\mu }^{\delta }}.
	\end{equation*}%
	From the embedding $ \mathcal{N}_{p,1,u}^{\sigma }\hookrightarrow \mathcal{N}_{p,\min(2,u),u}^{\sigma }\hookrightarrow \mathcal{E}_{p,2,u}^{\sigma } $ and Lemma \ref{Polya-Nikolskij1} the last term in \eqref{est-f1} can be
	estimated by%
	\begin{align*}
	c\sum_{j=J+1}^{\infty }2^{j\sigma }\big\|\mathcal{F}^{-1}\varphi _{j}\ast f%
	\big\|_{M_{u}^{p}} &\lesssim \sum_{j=J+1}^{\infty }2^{jn(\frac{1}{%
			q}-\frac{1}{p})+j\sigma }\big\|\mathcal{F}^{-1}\varphi _{j}\ast f\big\|_{%
		M_{v}^{q}} \\
	&\lesssim 2^{J(\frac{n}{q}-\frac{n}{p}+\sigma -s)}\big\|f\big\|_{\mathcal{N}%
		_{q,\infty ,v}^{s}},
	\end{align*}%
	since $s-\frac{n}{q}>\sigma -\frac{n}{p}$. Therefore,%
	\begin{align*}
	\big\|f\big\|_{\mathcal{E}_{p,2,u}^{\sigma }} &\leq c2^{J(\frac{n}{\delta }-%
		\frac{n}{p})+\sigma J}\big\|f\big\|_{M_{\mu }^{\delta }}+2^{J(%
		\frac{n}{q}-\frac{n}{p}+\sigma -s)}\big\|f\big\|_{\mathcal{N}_{q,\infty
			,v}^{s}} \\
	&=c2^{J(\frac{n}{\delta }-\frac{n}{p})+\sigma J}\left( \big\|f\big\|_{%
		M_{\mu }^{\delta }}+2^{J(\frac{n}{q}-\frac{n}{\delta }-s)}\big\|f%
	\big\|_{\mathcal{N}_{q,\infty ,v}^{s}}\right) ,
	\end{align*}%
	where the positive constant $c$\ is independent of $J$. We wish to choose $%
	J\in \mathbb{N}$ such that%
	\begin{equation*}
	\big\|f\big\|_{M_{\mu }^{\delta }}\approx 2^{J(\frac{n}{q}-\frac{n%
		}{\delta }-s)}\big\|f\big\|_{\mathcal{N}_{q,\infty ,v}^{s}},
	\end{equation*}%
	which is possible since $ \mathcal{N}_{q,\infty ,v}^{s}\hookrightarrow M_{\mu }^{\delta }$. Indeed, from Theorem \ref{Sobolev-embeddings} and \eqref{Morrey}, we get
	
	\begin{equation*}
	\mathcal{N}_{q,\infty ,v}^{s}\hookrightarrow\mathcal{E}_{q,\infty ,v}^{s} \hookrightarrow \mathcal{E}_{\delta,2 ,\mu}^{0}=M_{\mu }^{\delta },
	\end{equation*}
	becuase of $s-\frac{n}{q}>\sigma -\frac{n}{p}\geq -%
	\frac{n}{\delta }$. Thus%
	\begin{equation*}
	\big\|f\big\|_{\mathcal{E}_{p,2,u}^{\sigma }}\lesssim \big\|f\big\|_{%
		M_{\mu }^{\delta }}^{1-\theta }\big\|f\big\|_{\mathcal{N}%
		_{q,\infty ,v}^{s}}^{\theta }.
	\end{equation*}%
	Using \eqref{new-norm3}\ and \eqref{new-norm4} we arrive at
	the inequality%
	\begin{equation*}
	\big\|f\big\|_{\mathcal{\dot{E}}_{p,2,u}^{\sigma }}\lesssim \big\|f\big\|_{%
		M_{\mu }^{\delta }}^{1-\theta }\left( \big\|f\big\|_{M_{v}^{q}}+\big\|f\big\|_{\mathcal{\dot{N}}_{q,\infty ,v}^{s}}\right)
	^{\theta }.
	\end{equation*}%
	In this estimate replace $f$ by $f(\lambda \cdot )$ and using %
	\eqref{new-norm3} to obtain%
	\begin{equation*}
	\big\|f\big\|_{\mathcal{\dot{E}}_{p,2,u}^{\sigma }}\lesssim \big\|f\big\|_{%
		M_{\mu }^{\delta }}^{1-\theta }\left( \lambda ^{- s}\big\|f%
	\big\|_{M_{v}^{q}}+\big\|f\big\|_{\mathcal{\dot{N}}_{q,\infty
			,v}^{s}}\right) ^{\theta }.
	\end{equation*}%
	Taking $\lambda $ large enough we obtain \eqref{result-Morrey}-%
	\eqref{result-Morrey1}.
\end{proof}
\section*{\large Acknowledgments}We thank the referee for carefully reading the paper and for making several useful suggestions and comments.

 This work is found by the General Direction of Higher
	Education and Training under Grant No. C00L03UN280120220004 and by
	The General Directorate of Scientific Research and Technological Development,
	Algeria.

\end{document}